\documentclass[12pt,reqno]{amsart}
\usepackage{diagram}
\usepackage{amssymb}
\usepackage{graphicx}
\usepackage{amsmath}
\usepackage{a4wide}
\usepackage{version}
\usepackage{mathrsfs}
\usepackage{tikz}
\usepackage{mathtools}
\usetikzlibrary{arrows,decorations.markings, matrix}
\usepackage{tikz-cd}
\usepackage{hyperref}
\usepackage{float}
\usepackage[all]{xy}

\numberwithin{equation}{section}

\newtheorem{thm}{Theorem}[subsection]
\newtheorem{lem}[thm]{Lemma}
\newtheorem{prop}[thm]{Proposition}
\newtheorem{cor}[thm]{Corollary}

\theoremstyle{definition}
\newtheorem{definition}[thm]{Definition}
\newtheorem{example}[thm]{Example}
\newtheorem{rmk}[thm]{Remark}

\newcommand{\N}{\mathbb{N}}

\newcommand{\Z}{\mathbb{Z}}

\newcommand{\id}{\operatorname{id}}

\newcommand{\ov}{\overline}
\newcommand{\we}{\wedge}

\newcommand{\rk}{\operatorname{rk}}

\newcommand{\Hom}{\operatorname{Hom}}

\newcommand{\Ext}{\operatorname{Ext}}
\newcommand{\End}{\operatorname{End}}

\newcommand{\Aut}{\operatorname{Aut}}

\renewcommand{\a}{\alpha}
\renewcommand{\b}{\beta}
\newcommand{\om}{\omega}

\newcommand{\la}{\lambda}

\renewcommand{\P}{{\mathbb P}}
\newcommand{\C}{{\mathbb C}}

\newcommand{\wt}{\widetilde}
\newcommand{\ot}{\otimes}

\newcommand{\sub}{\subset}

\newcommand{\GL}{\operatorname{GL}}

\renewcommand{\mod}{\operatorname{mod}}

\newcommand{\und}{\underline}

\newcommand{\OO}{{\mathcal O}}

\newcommand{\G}{{\mathbb G}}

\newcommand{\hra}{\hookrightarrow}
\newcommand{\lan}{\langle}
\newcommand{\ran}{\rangle}

\renewcommand{\P}{{\mathbb P}}

\newcommand{\ga}{\gamma}
\newcommand{\de}{\delta}

\renewcommand{\ker}{\operatorname{ker}}
\newcommand{\im}{\operatorname{im}}

\newcommand{\fe}{{\mathbf e}}

\newcommand{\sspan}{\operatorname{span}}

\newcommand{\bd}{\mathbf{d}}
\newcommand{\lcm}{\mathrm{lcm}}
\newcommand{\bone}{\mathbf{1}}

\newcommand{\cU}{\mathcal{U}}
\newcommand{\cV}{\mathcal{V}}
\newcommand{\cN}{\mathcal{N}}

\newcommand{\cF}{\mathcal{F}}
\newcommand{\ol}{\overline}

\newcommand{\precdot}{\prec\mathrel{\mkern-5mu}\mathrel{\cdot}}

\newcommand{\cX}{\mathcal{X}}

\newcommand{\cok}{\mathrm{cok}}

\title{
Feigin-Odesskii brackets associated with Kodaira cycles and positroid varieties}
\author{Zheng Hua}
\author{Alexander Polishchuk}

\begin{document}

\begin{abstract}
We establish a link between open positroid varieties in the Grassmannians $G(k,n)$ and certain moduli spaces of complexes of vector bundles over Kodaira cycle $C^n$, using the shifted Poisson structure on the latter moduli spaces and relating them to a certain twist of the standard Poisson structure on $G(k,n)$. 
This link allows us to solve a classification problem for extensions of vector bundles over $C^n$. Based on this solution we further classify the symplectic leaves of all positroid varieties in $G(k,n)$ with respect to the twisted standard Poisson structure. Moreover, we get an explicit description of the moduli stack of symplectic leaves of $G(k,n)$ with the twisted standard Poisson structure as an open substack of the stack of  vector bundles on $C^n$.
\end{abstract}

\maketitle

\section{Introduction}
Positroid varieties (see Definition \ref{def: positroid}) are interesting subvarieties in complex Grassmannians 
that appeared in the study of total positivity and Poisson geometry. We refer to [KLS13] for
several other equivalent definitions.  They give a  stratification  
\begin{equation}\label{eq: strat}
G(k,n)=\bigsqcup_{f\in B(k,n)} X_f
\end{equation} where $X_f$ is the positroid variety with index $f$. 
The index set $B(k,n)$ is the set of bounded affine permutations (see definition in \cite{KLS} or Section \ref{sec: affper}). Unlike the Schubert cells, positroid varieties have very interesting topology \cite{GL-cat}. The cohomology classes of the closures of positroid varieties give the affine Stanley symmetric functions \cite{KLS}. Galashin and Lam \cite{GL-cluster} prove that the rings of functions on positroid varieties admit a cluster algebra structure in the sense of Fomin and Zelevinsky, which is closely related to the theory of total positivity developed by Lusztig and Postnikov. They also appear in mathematical physics \cite{Lam-physics} and algebraic combinatorics. 

Among various approaches toward positroid varieties, we are deeply influenced by the work of Goodearl and Yakimov. In \cite{GY}, Goodearl and Yakimov prove that the $T$-leaves of a partial flag variety of finite type with respect to the standard Poisson structure, where $T$ is the maximal torus of the corresponding complex reductive group, are the projected open Richardson varieties. In the special case of Grassmannians, the projected open Richardson varieties are precisely the positroid varieties. 

The work \cite{GY} suggests the question of constructing the ``moduli space" of symplectic leaves of $G(k,n)$.
In the case of the standard Poisson structure on a semi-simple complex Lie group the problem of constructing the ``space" of leaves was solved (only set theoretically)
by Kogan and Zelevinsky using ideas from the theory of cluster algebras \cite{KZ}.
 
In this paper, we consider
a {\it twist} of the standard Poisson structure on the Grassmannian $G(k,n)$ obtained by adding to the standard bivector a certain bivector coming from the torus action on $G(k,n)$.  We identify $G(k,n)$ equipped with the twisted standard Poisson structure with a certain moduli space of complexes of vector bundles on the Kodaira cycle $C^n$
by applying the machinery of shifted Poisson geometry developed in \cite{PTVV, CPTVV, MS} and our previous work \cite{HP-moduli, HP-matrix, HP-bih, HP-Bos}.

A connected projective curve $C$ is called Calabi-Yau (CY) if it is Gorenstein and $\omega_C\cong \OO_C$. Based on the seminal work of Pantev, Toen, Vaqui\'e and Vezzosi on shifted symplectic structure, we constructed a canonical $0$-shifted Poisson structure on the derived moduli stack of bounded complexes of vector bundles on $C$ up to chain isomorphisms, where we also can fix some terms and some cohomology sheaves (see \cite{HP-moduli, HP-Bos}). In the case when $C$ is an elliptic curve, it is known that this shifted Poisson structure is a derived lift of many classical Poisson structures of elliptic type including Feigin-Odesskii bracket on projective spaces \cite{FO}, Mukai-Bottacin Poisson structure on the regular part of moduli space of semi-stable sheaves on Poisson surfaces \cite{Bot}, etc. Derived algebraic geometric methods lead to many new results on these classical Poisson manifolds (see \cite{HP-moduli, HP-matrix, HP-bih, HP-Bos}). One notable component  $\cN_{L,k}$ of the moduli stack of complexes on $C$ parametrizes complexes given by an injective map $\OO_C^k\to V$ such that $V/\OO_C^k\cong L$, where $L$ is an ample line bundle on $C$ of degree $n\geq 1$, with some additional stability condition (see Definition \ref{def: NLk}). Regardless of the degeneration type of $C$, the coarse moduli space of $\cN_{L,k}$ is always the complex Grassmannian $G(k,n)$. However, the Poisson structure on $\cN_{L,k}$ depends on $C$. 

We prove that when $C=C^n$, the Kodaira cycle of $n$ smooth rational curves and $L$ is a line bundle which has degree one on each irreducible component, then the $0$-shifted Poisson structure on $\cN_{L,k}$ induces the twisted standard Poisson structure on $G(k,n)$ (see Theorem \ref{thm: comparison}). 


From our previous work, it is known that the sub moduli space of complexes $\OO_C^k\to V$ with fixed isomorphism type of $V$  is symplectic. It turns out that these sub moduli spaces are always connected (see Theorem \ref{thm: fibration}). Therefore the question of classifying symplectic leaves of the corresponding Poisson structure is equivalent to the classification problem of isomorphism types of vector bundles $V$ that are extensions of $L$ by $\OO_C^k$. We give a solution to this extension problem for $C=C^n$, based on a delicate correspondence between affine permutations and torus orbits on the stack of certain vector bundles (Theorem \ref{thm: Vtof}) and a new geometric intepretation of cyclic rank matrices of \cite{KLS}. We also determine the moduli space of symplectic leaves (see Theorem \ref{thm: main}). 
Namely, given $X_f$, we show that the moduli stack of symplectic leaves in $X_f$ is isomorphic to a torus orbit in the stack of vector bundles on $C^n$ (where the torus acts geometrically on $C^n$). Moreover, these tori glue into an open substack of the vector bundle stack (see Corollary \ref{cor: U++}). 
Actually, we also give a new proof of the theorem of Goodearl and Yakimov (\cite{GY}) classifying $T$-leaves in the case of Grassmannian, and of Yakimov's result about the dihedral group action on the set of positroid varieties (\cite[Corollary 1.2]{Ya}). In our geometric interpretation,
the dihedral group action
corresponds to the natural geometric action of this group on the Kodaira cycle $C^n$. This then leads 
to a Poisson (or anti-Poisson) automorphisms of the moduli stack of complexes (this follows from the general result
\cite[Theorem 7.5]{HP-Bos}).

The paper is organized as follows. In Section \ref{sec: FO} we discuss the construction of Feigin-Odesskii Poisson structures on projective spaces and its generalization to Grassmannians. In Section \ref{sec: standard} we prove our first result stating that the twisted standard Poisson structure on the Grassmannian $G(k,n)$ coincides with the Feigin-Odesskii Poisson structure on the component $\cN_{L,k}$, for $C=C^n$, the Kodaira cycle, and the line bundle $L$ on $C^n$ that has degree one on each irreducible component. In Section \ref{sec: vb} we recall basics on vector bundles on $C^n$ and begin the study of the stack $\cN_{L,k}$ and a related open substack $\cU^+_{L,k+1}$ of the stack of vector bundles on $C^n$ of rank $k+1$ and determinant $L$.
In particular, in Sec. \ref{geom-cyclic-sec}, we give a geometric definition of the cyclic rank matrix associated with a vector bundle. Section \ref{sec: affper} is devoted to various combinatorial results about affine permutations.   In Section \ref{sec: main} we establish the two main results of the paper. The first one is a correspondence between torus orbits 
on $\cU^+_{L,k+1}$ and certain affine permutations.  The second is the classification of extensions of $L$ by $\OO_C^k$ up to isomorphisms on $C=C^n$. In Sec. \ref{appl-sec} we show how this gives a classification theorem for symplectic leaves on positroid varieties.  

We work over the field of complex numbers.

\subsection*{Acknowledgment}
We thank Igor Burban, Jiang-hua Lu, Brent Pym, Xuhua He, Sasha Stolin and Raschid Abedin for many valuable discussions. 
Z.H. is partially supported by the GRF grant no. 17303420, 17305123 of University Grants Committee of Hong Kong SAR, China. 
A.P. is partially supported by the NSF grant DMS-2001224, Simons Travel Grant MP-TSM-00002745,
and within the framework of the HSE University Basic Research Program.

\section{Feigin-Odesskii Poisson structure}\label{sec: FO} 

\subsection{Bivector fields on the Grassmannian}\label{B-bivec-constr-sec}

Let $V$ be a finite dimensional vector space. A {\it quadratic bivector field} on $V$ is a linear map
$$b:{\bigwedge}^2 V^*\to S^2 V^*.$$
We can extend $b$ uniquely to a $\G_m$-invariant skew-symmetric bracket $\{\cdot,\cdot\}$ on the polynomial algebra $S(V^*)$. The corresponding bivector field restricts to $V\setminus \{0\}$ 
which descends to a global bivector field on $\P V$.

One way to construct a quadratic bivector field $b$ as above is by starting with a polylinear map
$$B(\cdot,\cdot;\cdot,\cdot):V^*\ot V\ot V^*\ot V\to \C,$$
with the following skew-symmetry:
\begin{equation}\label{B-skew-eq}
B(v_2^*,v_2;v_1^*,v_1)=-B(v_1^*,v_1;v_2^*,v_2),
\end{equation}
and then define $b(v_1^*,v_2^*)$ to be the quadratic form $B(v_1^*,x,v_2^*,x)$. 

There is a similar construction giving global bivector fields on Grassmannians.
Namely, let $W$ be another vector space of dimension $d$. 
Then we can identify the Grassmannian $G(d,V)$ with the quotient $\GL(W)\backslash X $, where
$$X=X_{W,V}\sub\Hom(W,V)$$ 
is the space of injective maps.

Let $p:X\to G(d,V)$ denote the natural projection. We have a canonical map
$${\bigwedge}^2T_X\to p^*{\bigwedge}^2T_{G(d,V)}$$
which induces a map
$$H^0(X,{\bigwedge}^2T_X)^{\GL(W)}\to H^0(X,p^*{\bigwedge}^2T_{G(d,V)})^{\GL(W)}\simeq H^0(G(d,V),{\bigwedge}^2T_{G(d,V)}).$$
Thus, starting with a $\GL(W)$-invariant quadratic bivector field $b_W$ on $\Hom(W,V)$, we get a $GL(W)$-invariant bivector field on $X$,
which induces a bivector field $\ov{b}_W$ on $\GL(W)\backslash X \simeq G(d,V)$.

Now we observe that starting with a polylinear map $B$ as above, satisfying the skew-symmetry
\eqref{B-skew-eq}, we can construct a $\GL(W)$-invariant quadratic bivector field $b_W$ on $\Hom(W,V)=W^*\ot V$ by setting
$$b_W(w_1\ot v_1^*,w_2\ot v_2^*)(\phi):=B(v_1^*,\phi(w_2),v_2^*,\phi(w_1)),$$
where $\phi\in \Hom(W,V)$, $w_i\in W$, $v_i^*\in V^*$.
The corresponding bivector field at the point $\phi(W)$ of the Grassmannian $G(d,V)$, is the restriction of the skew-symmetric form 
$b_W(\cdot,\cdot)(\phi)$ to
$W\ot \phi(W)^\perp$, i.e., we only consider $v_1^*,v_2^*\in \phi(W)^\perp$. Here $\phi(W)^\perp:=\{v^*\in V| v^*(\phi(w))=0 \text{ for any $w\in W$}\}$.

\subsection{Feigin-Odesskii brackets on the Grassmannians}

Let $\xi$ be an endosimple vector bundle on a CY curve $C$. Let us fix a trivialization $\OO_C\simeq \om_C$.
Then we have the corresponding FO bracket $\Pi_\xi$ on $\P V$, where
$V=\Ext^1(\xi,\OO_C)\simeq H^0(C,\xi)^*$. Let $\cN_{\xi,1}$ be the moduli stack of complexes of vector bundles $\OO_C\to E$ such that $E/\OO_C\cong \xi$ and $E\not\cong \xi\oplus \OO_C$.  The coarse moduli space of $\cN_{\xi,1}$ is isomorphic to $\P V$. The Poisson structure $\Pi_\xi$ is the classical shadow of the $0$-shifted Poisson structure on the moduli stack of complexes constructed and studied in \cite{HP-moduli, HP-matrix, HP-bih, HP-Bos}. 

More generally, for a vector space $W$, we have the FO bracket $\Pi_{W,\xi}$ on the space of extensions
\begin{equation}\label{WV-ext}
0\to W^*\ot \OO_C\to E\to \xi\to 0,
\end{equation}
corresponding to injective maps $W\to V=\Ext^1(\xi,\OO_C)$, modulo the $\GL(W)$-action, which can be identified with $G(d,V)$, where $d=\dim W$. Similarly, $G(d,V)$ is the coarse moduli space of a stack $\cN_{\xi,k}$ (see Definition \ref{def: NLk} below) and $\Pi_{W,\xi}$ is the classical shadow of a $0$-shifted Poisson structure. 
We will need a formula for $\Pi_{W,\xi}$ in terms of a certain triple Massey product.

\begin{lem}\label{FO-bracket-lem} 
For an extension corresponding to $\phi\in W^*\ot H^1(\xi^{-1})$, we identify the cotangent space $T^*_\phi G(d,V)$ with $W\ot \phi(W)^\perp\sub W\ot H^0(\xi)$,
where we view $\phi$ as a linear map $W\to H^1(\xi^{-1})$. In terms of the Serre duality between $H^0(\xi)$ and $H^1(\xi^{-1})$, the skew-symmetric pairing on $T^*_\phi G(d,V)$ associated with $\Pi_{W,\xi}$ is given by
$$\Pi_{W,\xi}(s_1\we s_2)=\pm \lan MP(s_1,\phi,s_2),\phi\ran,$$
where $s_1,s_2\in W\ot \phi(W)^\perp$, and we use the triple Massey product
$$MP: \Hom(W^\vee\ot \OO_C,\xi)\ot\Ext^1(\xi,W^\vee\ot \OO_C)\ot\Hom(W^\vee\ot\OO_C,\xi)\to \Hom(W^\vee\ot \OO_C,\xi).$$
\end{lem}

\begin{proof}
The proof is a straightforward generalization of \cite[Lem.\ 2.1]{HP-bih}. Let us recall the definition of the FO bracket in this case.
Given an extension \eqref{WV-ext}, let $\und{\End}(E,W^\vee\ot \OO_C)$ denote the bundle of endomorphisms of $E$ preserving $W^\vee\ot \OO_C$.
It fits into an exact sequence
\begin{equation}\label{rel-End-seq1}
0\to W^\vee\ot \xi^\vee\to \und{\End}(E,W^\vee\ot \OO_C)\to \und{\End}(\xi)\oplus \End(W)\ot\OO_C\to 0.
\end{equation}

\noindent
{\bf Step 1}. We claim that the induced map (obtained by applying $\Hom(\cdot,\OO_C)$ to \eqref{rel-End-seq1}),
$\Hom(\und{\End}(E,W^\vee\ot \OO_C),\OO_C)\to W\ot H^0(\xi)$, has $W\ot \phi(W)^\perp$ as an image.
Indeed, the exact sequence of $\Ext^*(\cdot,\OO_C)$ applied to \eqref{rel-End-seq1} identifies this image with 
the kernel of the connecting map
$$W\ot H^0(\xi)=\Hom(W^\vee\ot \xi^\vee,\OO_C)\to \Ext^1(\und{\End}(\xi)\oplus \End(W)\ot\OO_C,\OO_C).$$
By Serre duality this map is dual to the connecting homomorphism
\begin{equation}\label{End-xi-connecting-hom-eq}
H^0(\und{\End}(\xi)\oplus \End(W)\ot\OO_C)\to W^\vee\ot H^1(\xi^\vee)
\end{equation}
from the long sequence of cohomology associated with \eqref{rel-End-seq1}.
It remains to show that the image of the latter connecting homomorphism is $W^\vee\ot \phi(W)$.

First, we observe that the restriction of \eqref{End-xi-connecting-hom-eq} to $H^0(\End(W)\ot \OO_C)$ is
the connecting homomorphism associated with the exact sequence
$$0\to W^\vee\ot \xi^\vee\to W^\vee\ot E^\vee\to \End(W)\ot \OO_C\to 0$$
obtained by tensoring the dual of \eqref{WV-ext} with $W^\vee$. Hence, the image of this restriction is $W^\vee\ot \phi(W)$.

Next, we claim that the summand $H^0(\und{\End}(\xi))\simeq H^0(\OO_C)$ does not increase the image of \eqref{End-xi-connecting-hom-eq}.
Indeed, this immediately follows from the fact that the extension \eqref{rel-End-seq1} splits over the diagonal
embedding $\OO_C\to \und{\End}(\xi)\oplus \End(W)\ot \OO_C$. 

\noindent
{\bf Step 2}. The Poisson bivector field $\Pi_{W,\xi}|_\phi$ is uniquely determined by the commutative diagram
\begin{diagram}
\Hom(\und{\End}(E,W^\vee\ot \OO_C),\OO_C)&\rTo{\de}& W^\vee \ot H^1(\xi^\vee)\\
\dTo{\a_1}&&\dTo{\a_2}\\
W\ot \phi(W)^\perp&\rTo{\Pi_{W,\xi}|_\phi}& W^\vee \ot H^1(\xi^\vee)/\im(\phi)
\end{diagram}
where $\a_1$ is a surjective map from Step 1, and $\de$ is obtained by applying $R\Hom(\cdot,\OO_C)$ to the sequence
\begin{equation}\label{rel-End-seq2}
0\to \und{\End}(E,W^\vee\ot \OO_C)\to\und{\End}(E)\to W\ot \xi\to 0.
\end{equation}

Indeed, recall that $\Pi_{W,\xi}$ is induced by the Poisson bivector field on the moduli of complexes, which has $H^1(\und{\End}(E,W^\vee\ot \OO_C))$ as the tangent space,
and the latter is given by the composition
$$H^0(\und{\End}(E,W^\vee\ot \OO_C)^\vee)\rTo{\de} H^1(W^\vee \ot \xi^\vee)\to H^1(\und{\End}(E,W^\vee\ot \OO_C))$$
where the second arrow comes from the long exact sequence associated with \eqref{rel-End-seq1}. Now the assertion follows from the fact
that this second arrow is exactly the composition of the projection $\a_2$ with the embedding of tangent spaces
$$W^\vee \ot H^1(\xi^\vee)/\im(\phi)\hra H^1(\und{\End}(E,W^\vee\ot \OO_C),$$
together with the fact that $\a_1$ is identified with the dual of this embedding.

\noindent
{\bf Step 3}. The morphism of exact sequences
\begin{diagram}
0 &\rTo{}& \xi^\vee\ot E&\rTo{}&\und{\End}(E)&\rTo{}& W\ot E &\rTo{}& 0\\
&&\dTo{}&&\dTo{}&&\dTo{}\\
0 &\rTo{}&\und{\End}(E,W^\vee\ot\OO_C)&\rTo{}&\und{\End}(E)&\rTo{}&W\ot \xi&\rTo{}&0
\end{diagram}
leads to a commutative square
\begin{diagram}
\Hom(\und{\End}(E,W^\vee\ot \OO_C),\OO_C)&\rTo{\de}& W^\vee \ot H^1(\xi^\vee)\\
\dTo{\b_1}&&\dTo{\b_2}\\
H^0(\xi\ot E^\vee)&\rTo{\de'}& W^\vee \ot H^1(E^\vee)
\end{diagram}
where $\de'$ is obtained by applying $R\Hom(\cdot,\OO_C)$ to the first row.

\noindent
{\bf Step 4}. The Massey product map $MP(\phi,\cdot,\phi)$ given by the composable arrows
$$\xi\rTo{[1]} W^\vee \ot \OO_C\to \xi \rTo{[1]} W^\vee\ot \OO_C$$
fits into a commutative square
\begin{diagram}
H^0(\xi\ot E^\vee)&\rTo{\de'}& W^\vee \ot H^1(E^\vee)\\
\dTo{\ga_1}&&\dTo{\ga_2}\\
W\ot \phi(W)^\perp&\rTo{MP(\phi,\cdot,\phi)}& W^\vee \ot H^1(\xi^\vee)/\im(\phi)
\end{diagram}

\noindent
{\bf Step 5}. Next, we observe that $\a_1=\b_1\ga_1$ and $\a_2=\b_2\ga_2$,
which identifies the map $\Pi_{W,\xi}|_{\phi}$ with $MP(\phi,\cdot,\phi)$.
Finally, we use cyclic symmetry to relate this to $MP(\cdot,\phi,\cdot)$, exactly as in the proof of \cite[Lem.\ 2.1]{HP-bih}.
\end{proof}

\begin{prop}\label{FO-grassm-prop} 
The FO bracket $\Pi_{W,\xi}$ on the Grassmannian $\GL(W)\backslash X \simeq G(d,V)$
is obtained from the polylinear map
$$V^*\ot V\ot V^*\ot V\to \C: s_1\ot e\ot s_2\ot e'\mapsto \lan m_3(s_1,e,s_2),e'\ran,$$
where we identify $V^*$ with $H^0(\xi)$,
by the construction of Sec.\ \ref{B-bivec-constr-sec}.
\end{prop}

\begin{proof} This follows from the formula of Lemma \ref{FO-bracket-lem}.
Namely, let $(e_i)$ be a basis of $W$, and let 
$$\phi=\sum_i e_i^*\ot v_i, \ s_1=e_i\ot \a, \ s_2=e_{i'}\ot \b,$$
where $v_i\in H^1(\xi^{-1})$, $\a,\b\in \lan v_1,\ldots,v_d\ran^\perp\sub H^0(\xi)$.
Then we have
$$MP(s_1,\phi,s_2)=e_i\ot MP(\a,v_{i'},\b)\in W\ot H^0(\xi).$$
Hence,
$$\Pi_{W,\xi}(s_1\we s_2)=\lan MP(\a,v_{i'},\b), v_i\ran=\lan m_3(\a,v_{i'},\b),v_i\ran.$$
\end{proof}

\section{A modular construction of Standard Poisson bracket}\label{sec: standard}

\subsection{The Massey products for the Kodaira cycles}

Let $C=\bigcup C_i$ denote the Kodaira cycle of length $n$, where $C_i=\P^1$, $C_i$ is glued with $C_{i+1}$ at the points $\infty\in C_i$ and $0\in C_{i+1}$.
Let $p_i\in C_i\setminus\{0,\infty\}$ be the point with the coordinate $z_i=1$. We set 
$$L:=\OO_C\Big(\sum_{i=1}^n p_i\Big).$$
We will compute the cohomology of $L$, $L^{-1}$ and $\OO_C$ using Cech resolutions with respect to the covering $(U_0,U_1)$, where
$U_0$ is the formal neighborhood of $\{p_1,\ldots,p_n\}$, $U_1=C\setminus\{p_1,\ldots,p_n\}$.
We will use $x_i=z_i-1$ as a formal parameter near $p_i$.
Then the relevant Cech complexes have form
$$C(\OO_C): \bigoplus_i \C[\![x_i]\!] \oplus \OO_C(U_1)\to \bigoplus_i \C(\!(x_i)\!),$$
$$C(\xi): \bigoplus_i \C[\![x_i]\!]\cdot x_i^{-1} \oplus \OO_C(U_1)\to \bigoplus_i \C(\!(x_i)\!),$$
$$C(\xi^{-1}): \bigoplus_i \C[\![x_i]\!]\cdot x_i \oplus \OO_C(U_1)\to \bigoplus_i \C(\!(x_i)\!),$$
where the Cech differential has form $\de(a_0,a_1)=a_1|_{U_{01}}-a_0|_{U_{01}}$.

We denote by $e_i$ the natural idempotents in each direct sum above.
Let us set $f_i[n]=(1+x_i^{-1})x_i^{-n}\cdot e_i$. Then for $n\ge 1$, $f_i[n]$ is a well defined element of $\OO_C(U_1)$.
We also have elements 
$$s_i:=(1+x_{i-1}^{-1})e_{i-1}-x_i^{-1}e_i\in \OO_C(U_1),$$
$i\in \Z/n\Z$, which extend to global sections of $L$ and form a basis of $H^0(L)$.
We have a unique isomorphism $\tau:H^1(\OO_C)\to \C$ sending the class of $x_i^{-1}e_i\in C^1(\OO_C)$ to $1$, for each $i$.
Also, the classes $x_i^{-1}e_i\in C^1(L^{-1})$ form a basis of $H^1(L^{-1})$ which is dual to the basis $(s_i)$ of $H^0(L)$ with
respect to the perfect pairing 
$$H^0(L)\ot H^1(L^{-1})\to \C: \lan \a,\b\ran=\tau(\a\b).$$
Indeed, this follows from the relations
\begin{equation}\label{de-fi-eqs}
\begin{array}{l}
x_i^{-1}e_i\cdot s_i=-x_i^{-2}e_i=x_i^{-1}e_i-\de(0,f_i[1]),\\
x_i^{-1}e_i\cdot s_{i+1}=\de(0,f_i[1]),\\
x_i^{-1}e_i=x_n^{-1}e_n+\de(\sum_{i\le j<n}e_j,\sum_{i<j\le n}s_j).
\end{array}
\end{equation}

Now assume we are given elements 
$$\a=\sum_i a_i s_i\in H^0(L), \ \ \b=\sum_i b_is_i\in H^0(L), e=\sum_i \la_i x_i^{-1}e_i\in H^1(L^{-1})$$
such that 
$$\lan \a,e\ran=\sum a_i\la_i=0, \ \ \lan \b,e \ran=\sum b_i\la_i=0.$$
Then we can consider the triple Massey product $MP(\a,e,\b)\in H^0(\xi)/\sspan(\a,\b)$. Hence,
if we have another element $e'=\sum_i \mu_i x_i^{-1}e_i\in H^1(L^{-1})$ such that 
$\lan \a,e'\ran=\lan \b,e'\ran=0$, then $\lan MP(\a,e,\b),e'\ran$ is well defined.

\begin{lem}\label{MP-Kod-lem}
One has
$$\lan MP(\a,e,\b),e'\ran=\sum_{i<j} (a_ib_j-b_ia_j)\la_i\mu_j.$$
\end{lem}

\begin{proof}
We use the natural representative of $e$ in $C^1(L^{-1})$.
By definition, 
$$MP(\a,e,\b)=\de^{-1}(\a\cdot e) \cdot \b-\de^{-1}(\b\cdot e)\cdot \a.$$
We have
\begin{align*}
&\a\cdot e=-\sum_i  a_i\la_i x_i^{-2}e_i + \sum_i a_i\la_{i-1}(x_{i-1}^{-1}+x_{i-1}^{-2})e_{i-1}=\\
&\sum_i a_i\la_i x_i^{-1}e_i+\de(0,-\sum_i a_i\la_i f_i[1]+\sum_i a_i\la_{i-1}f_{i-1}[1])=\\
&\de\Big(\sum_{1\le i\le j<n}a_i\la_ie_j,\sum_{1\le i<j\le n}a_i\la_i s_j-\sum_i a_i\la_i f_i[1]+\sum_i a_i\la_{i-1}f_{i-1}[1]\Big),
\end{align*}
where we used relations \eqref{de-fi-eqs}.
Hence,
$$\de^{-1}(\a\cdot)\cdot \b=(\sum_{1\le i\le j<n}a_i\la_ie_j)\cdot(\sum_{i=1}^n b_i s_i)=
\sum_{j=1}^{n-1}e_j\cdot (\sum_{i=1}^ja_i\la_i)(-b_jx_j^{-1}+b_{j+1}(1+x_j^{-1})),$$
and so
$$MP(\a,e,\b)=\sum_{j=1}^{n-1} e_j\cdot [(\sum_{i=1}^ja_i\la_i)(-b_jx_j^{-1}+b_{j+1}(1+x_j^{-1})-
(\sum_{i=1}^jb_i\la_i)(-a_jx_j^{-1}+a_{j+1}(1+x_j^{-1}))]$$
in $H^0(\xi)/\sspan(\a,\b)$.
Now it is easy to check that 
$$MP(\a,e,\b)=\sum_{j=2}^n s_j\cdot \Big[b_j\cdot (\sum_{i=1}^{j-1} a_i\la_i)-a_j\cdot (\sum_{i=1}^{j-1} b_i\la_i)\Big],$$
where the term with $s_n$ can be omitted (since the corresponding coefficient is zero).
This immediately leads to the claimed formula.
\end{proof}

The automorphism group of $C$, denoted by $\Aut(C)$ contains a torus $\ol{T}$ of dimension $n$. Let $t=(t_1,\ldots,t_n)\in \ol{T}$. It acts on the local coordinate $z_i$ by $t\cdot z_i=t_i z_i$. Let $w_i$ be a regular point on $C_i$ for $i=1,\ldots,n$, therefore correspond to unique complex numbers under the local coordinate $z_i$. The line bundle $L=\OO_C(w_1+\ldots+w_n)$ has degree one on each irreducible component. One can show that any line bundle of degree one on each component can be expressed in this form. Then
\[
t^*L=\OO_C(\sum_{i=1}^n t_i^{-1}w_i).
\] We claim that $\OO_C(w_1+\ldots+w_n)\cong \OO_C(w^\prime_1+\ldots+w^\prime_n)$ if and only $\prod_{i=1}^n w_i=\prod_{i=1}^n w^\prime_i$.  This can be proved as follows. On the irreducible component $C_i$, we take the rational function $\frac{z_i-w_i}{z_i-w^\prime_i}$ which induces an isomorphism $\OO_{C_i}(w_i)\cong \OO_{C_i}(w^\prime_i)$. Since we glue $\infty\in C_i$ with $0\in C_{i+1}$, the collection $\{\frac{z_i-w_i}{z_i-w^\prime_i}|i=1,\ldots,n\}$ of rational functions glue to a global rational function on $C$ if and only if the equality of product holds. Now we take $w^\prime_i=t_i^{-1}w_i$. Then $\prod_{i=1}^n w_i=\prod_{i=1}^n w^\prime_i$ when $\prod_{i=1}^nt_i=1$. Denote by $T\subset \ol{T}$ the sub torus satisfying $\prod_{i=1}^nt_i=1$. Suppose $t\in T$, we may use the above rational function to construct an isomorphism between $t^*L$ and $L$. Clearly the isomorphism depends on $t$ holomorphically. 

\begin{definition}\label{def: NLk}
Let $C$ be the Kodaira cycle with $n$ irreducible components and $L$ be a line bundle on $C$ with degree one on each component. Fix integer $0<k<n$ and a Noetherian scheme $S$. Define $\cN_{L,k}$ to be the moduli stack whose $S$-points are  complexes of vector bundles $\OO_{C\times S}^k\to \cV$ such that $\cV$ is a vector bundle on $C\times S$ and for $s\in S$
\begin{enumerate}
\item[$(1)$] $\cV_s/\OO_C^k\cong L$;
\item[$(2)$] the connecting morphism $H^0(L)\to H^1(\OO_C^k)$ is surjective.
\end{enumerate}
\end{definition}
The torus $T$ acts on $\cN_{L,k}$ as follows. Given a short exact sequence
\[\xymatrix{
\OO_C^k\ar[r]^a & V\ar[r]^b & L }
\] that represents a $\C$-point of $\cN_{L,k}$ and $t\in T$, we get an exact sequence
\[\xymatrix{
t^*\OO_C^k\ar[r]^{t^*a} & t^*V\ar[r]^{t^*b} & t^*L}.
\] Since $t^*\OO_C\cong \OO_C$ and $t^*L\cong L$, we obtain a complex $\OO_C^k\to t^*V$ with $t^*V/\OO_C\cong L$. Condition $(2)$ is preserved under $t^*$ since it is an equivalence of categories. 

\subsection{Twisted standard Poisson bracket on the Grassmannian}

Let $V$ be a finite-dimensional vector space, $W$ a vector space of dimension $d$, $X=X_{W,V}\sub \Hom(W,V)$ the space of injective maps, 
$G(d,V)=GL(W)\backslash X_{W,V}$ the corresponding Grassmannian.

The natural action of $\GL(V)$ on the Grassmannian $G(d,V)$ induces a linear map
$$\chi:\End(V)/\C\cdot \id\to H^0(G(d,V),T_{G(d,V)}).$$
Let us choose a basis $(e_i)_{1\le i\le n}$ of $V$, and let $E_{ij}\in \End(V)$ denote the corresponding elementary matrices.
The standard Poisson bivector field on $G(d,V)$ is given by the bivector field 
$$\pi_{d,V}=\sum_{1\le i<j\le n}\chi(E_{ij})\we \chi(E_{ji}).$$

We define the {\it twisted standard Poisson bracket on the Grassmannian} by the bivector
\begin{equation}\label{mod-st-Pois-br}
\pi'_{d,V}=\sum_{1\le i<j\le n}\chi(E_{ij})\we \chi(E_{ji})+\sum_{1\le i<j\le n}\chi(E_{ii})\we \chi(E_{jj}).
\end{equation}
\begin{rmk}
Let $(X,\pi)$ be a smooth Poisson variety such that $\pi$ is invariant under an action of a commutative algebraic group $G$ on $X$. The natural map from the Lie algebra $\rm{Lie}(G)$ to vector fields on $X$, induces a map 
$${\bigwedge}^2{\rm Lie}(G)\to H^0(X,{\bigwedge}^2 T_X): \theta\mapsto \pi_\theta.$$
For any $\theta\in {\bigwedge}^2{\rm Lie}(G)$, the bivector field $\pi+\pi_\theta$ is still Poisson. We call $\pi+\pi_\theta$ the $\theta$-twist of $\pi$. The $\theta$-twist does not change the $G$-leaves. In the above example, we take $(X,\pi)$ to be $(G(d,V), \pi_{d,V})$ and $G=T$ to be the maximal torus of $GL(V)$.
\end{rmk}

\begin{thm}\label{thm: comparison}
The following Poisson bivector fields on $G(d,V)$ are the same:
\begin{itemize}
\item the twisted standard bracket $\pi'_{d,V}$;
\item the bracket induced via the construction of Sec.\ \ref{B-bivec-constr-sec} by the polylinear map
$$B'_{st}:V^*\ot V\ot V^*\ot V\to \C$$
given by
$$B'_{st}(e_i^*\ot e_j\ot e_k^*\ot e_l)=\begin{cases}\de_{ij}\de_{kl}-\de_{il}\de_{jk}, & j<l, \\ -\de_{ij}\de_{kl}+\de_{il}\de_{jk}, & j>l \\  0, & j=l;\end{cases}$$
\item the rescaled FO bracket $2\Pi_{W,L}$ associated with  any $W$ and a line bundle $L$ on the Kodaira cycle of length $n$, where $L$ has degree $1$ on each component.
\end{itemize}
\end{thm}

\begin{proof}

For $A,B\in \End(V)$, the bivector field $\chi(A)\wedge \chi(B)$ at the point of the Grassmannian $W\sub V$ is given by the pairing
$$\chi(A)_W\we \chi(B)_W (w_1\ot v_1^*,w_2\ot v_2^*)=\lan v_1^*,A(w_1)\ran \lan v_2^*,B(w_2)\ran-\lan v_1^*,B(w_1)\ran \lan v_2^*,A(w_2)\ran.$$

In other words, $\chi(A)\we \chi(B)$ is induced by the polylinear form
$$B_{A,B}(v_1^*,x,v_2^*,y)=\lan v_1^*,A(y)\ran \lan v_2^*,B(x)\ran-\lan v_1^*,B(y)\ran \lan v_2^*,A(x)\ran.$$

Thus, the first sum in the definition of $\pi'_{d,V}$ corresponds to the form
\begin{align*}
\sum_{i'<j'} B_{E_{i'j'},E_{j'i'}}(e_i^*\ot e_j\ot e_k^*\ot e_l)
=\sum_{i'<j'}
[\lan e_i^*,E_{i'j'}(e_l)\ran \lan e_k^*,E_{j'i'}(e_j)\ran - \lan e_i^*,E_{j'i'}(e_l)\ran \lan e_k^*,E_{i'j'}(e_j)\ran].
\end{align*}
The first term is nonzero only if $j<l$ and is equal to $\de_{ij}\de_{kl}$. The second term is nonzero only if $j>l$ and is equal to $-\de_{ij}\de_{kl}$.
Similarly, the second sum in the definition $\pi'_{d,V}$ gives
\begin{align*}
\sum_{i'<j'} B_{E_{i'i'},E_{j'j'}}(e_i^*\ot e_j\ot e_k^*\ot e_l)
=\sum_{i'<j'}
[\lan e_i^*,E_{i'i'}(e_l)\ran \lan e_k^*,E_{j'j'}(e_j)\ran - \lan e_i^*,E_{j'j'}(e_l)\ran \lan e_k^*,E_{i'i'}(e_j)\ran].
\end{align*}
The first term is nonzero only if $j>l$ and is equal to $\de_{il}\de_{jk}$. The second term is nonzero only if $j<l$ and is equal to $-\de_{il}\de_{jk}$.
This proves that $\pi'_{d,V}$ is induced by $B'_{st}$.

The fact that the multlinear form $B'_{st}$ induces the rescaled FO bracket $2\Pi_{W,L}$ follows from Proposition \ref{FO-grassm-prop} and Lemma \ref{MP-Kod-lem}.

\end{proof}

\section{Vector bundles over Kodaira cycles}\label{sec: vb} 
\subsection{Classification of indecomposable bundles}

Recall that the Kodaira cycle with $n$ irreducible components is denoted by $C^n=\bigcup_{i\in\Z/n\Z} C^n_i$  where $C^n_i\cong \P^1$ and $C^n_i$ intersects transversely with $C^n_{i+1}$.  When the number $n$ is fixed in the context, we will omit the superscript $n$ from the notation. Let $V$ be a  vector bundle on $C$ of rank $r$. Then
\[
V|_{C_i}=\bigoplus_{j=0}^{r-1} \OO_C(d_{i+jn}).
\]
In particular, if $V$ is a line bundle then we say it has degree vector $\bd=(d_1,\ldots,d_n)$.  Given a degree vector $\bd$, line bundles with degree vector $\bd$ is classified by $\C^*$. For $\lambda\in \C^*$, there is a unique line bundle $L(\bd,\lambda)$. In fact one can show that $L(\bd,\lambda)\cong \OO_C(D)$ for some effective divisor $D=D_1+\ldots+D_n$ such that $D_i$ is supported on a collection of regular points $w^i_1,\ldots,w^i_{d_i}$ (of $C$) on $C_i$, and we have
\[
\lambda=\prod_{i=1}^n\prod_{j=1}^{d_i} w^i_{j}.
\]
Given a sequence of integers
\[
\bd=(d_1,\ldots,d_{rn})
\] and $\wt{\lambda}\in \C^*$, let $\wt{L}(\bd,\wt{\lambda})$ be the  line bundle on $C^{rn}$ such that $\wt{L}(\bd,\wt{\lambda})|_{C^{rn}_i}=\OO_{\P^1}(d_i)$. The cyclic group $\mu_{rn}:=\{\exp(\frac{j}{rn}2\pi i)|j\in\Z\}$ acts on $C^{rn}$ by cyclic rotation. Denote by $\pi_r: C^{rn}\to C^n$ the $r$ to $1$ cyclic cover by taking the quotient by $\mu_r:=\{\exp(\frac{2\pi i}{r})|i\in\Z\}$  and set $B(\bd,\wt{\lambda}):=(\pi_r)_*\wt{L}(\bd,\wt{\lambda})$. For $m\in\N$, there is a unique rank $m$ indecomposable bundle $\cF_m$ that is a $m$-th iterated extension of $\OO_C$. We write
\[
B(\bd,m,\wt{\lambda}):=B(\bd,\wt{\lambda})\ot\cF_m.
\] It is a vector bundle of rank $rm$ and degree $m\cdot \sum_{i=1}^{rn} d_i$.  
\begin{thm}\cite[Theorem 2.12]{DG}\cite[Theorem 3.15]{BBDG}
If $\bd$ is non-periodic, i.e. not of the form $\fe^s=\fe\ldots\fe$ where $\fe$ is an integral sequence of length $\frac{rn}{s}<rn$ and $s$ divides $r$,  then $B(\bd, m,\wt{\lambda})$ is indecomposable. Moreover $B(\bd, m,\wt{\lambda})\cong B(\bd', m',\wt{\lambda}')$ if and only if
\begin{equation}\label{eq: d'}
\bd'=(d_{n+1},\ldots,d_{2n},\ldots,d_{rn},d_1,\ldots,d_n)
\end{equation} and $m=m', \wt{\lambda}=\wt{\lambda}'$. Let $V$ be an indecomposable vector bundle on $C=C^n$. Then there exists $\bd, \wt{\lambda}$ and $m\in\N$ such that  $V\cong B(\bd, m, \wt{\lambda})$.
\end{thm}
We will call $m$ the \emph{multiplicity} of the indecomposable $B(\bd,m,\lambda)$.
Set $\tau\bd:=\bd'$ in equation \ref{eq: d'}. Consider $r\times rn$ matrix
\[
A(\bd,r)=\left(\begin{array}{cccc}
\bd \\ \tau\bd\\ \ldots \\ \tau^{r-1}\bd
\end{array}\right)
\]  up to permutation of rows. The above theorem shows that there is a unique equivalence class of matrices $A(\bd,r)$ (up to row permutations) associated to an indecomposable bundle $V=B(\bd,\wt{\lambda})$.  And the equivalence class only depends on $\bd$. Sometimes we also write $A(V)$ by an abuse of notation.

Suppose now $V=V_1\oplus \ldots\oplus V_p$ such that $V_i$ are pairwise non isomorphic indecomposables, of multiplicity $m_i=1$ and rank $r_i$ for $i=1,\ldots,p$. Let $h=\lcm(r_1,\ldots,r_p)$ and $h=s_ir_i$ for $i=1,\ldots,p$. Consider an $r\times hn$ matrix
\[
A(V):=\left(
\begin{array}{cccc}
A(V_1) & \ldots & A(V_1)\\
A(V_2) & \ldots & A(V_2)\\
&\ldots &\\
A(V_p) & \ldots & A(V_p)
\end{array}
\right)
\] where the $i$-th row has $s_i$ blocks. Similarly to the indecomposable case, we consider  $A(V)$ up to row permutations. In the later applications, it is more convenient to consider an $r\times \infty$ matrix that is the periodic extension of $A(V)$ in both left and right directions. By an abuse of notations, we use the same symbol $A(V)$ to represent this infinite matrix.

\subsection{Substack $\cU^+_{L, k+1}$}
Denote by  $\bone$ the vector $(1,\ldots,1)$. 
\begin{lem}
 Suppose $B(\bd,m, \wt{\lambda})$ is a summand of an extension bundle $V$ of  $L(\bone,\lambda)$ by $\OO_C^{r-1}$. Then $d_i=0,1$ for $i=1,\ldots,rn$. If in addition $H^1(V)=0$ then $m=1$ and $\bd$ is nonzero.
\end{lem}
\begin{proof}
 Given an exact sequence $\OO_C^{r-1}\to V\to L(\bone,\lambda)$, let $B:=B(\bd,m,\wt{\lambda})$ be a summand of $V$. By restricting the exact sequence to $C_i$, we get a $V|_{C_i}\cong \OO^{r-1}\oplus \OO(1)$. So $0\leq d_i\leq 1$ for $i=1,\ldots,rn$.  The assumption $H^1(V)=0$ implies that $H^1(B)=0$. Suppose $B=B(\bd,m,\wt{\lambda})$ then either $\bd$ is nonzero or $\bd=0$ and $\wt{\lambda}=1$. The latter case is impossible since then $V/\OO_C^{r-1}$ has $B$ as a summand. In the first case we conclude that $m=1$ by the determinant condition.
\end{proof}
Fix $L=L(\bone,\lambda)$ for some fixed $\lambda\in\C^\times$.
Let $S$ be a Noetherian scheme. We denote by $\cU^+_{L,k+1}$ the stack whose $S$ points are vector bundles $\cV$ of rank $k+1$  on $C\times S$ such that for $s\in S$
\begin{enumerate}
\item[$(1)$] $\det\cV_s$ is isomorphic to $L$;
\item[$(2)$]  for every indecomposable summand of $\cV_s$, the degree vector $\bd$ has entries $0,1$ and is a nonzero vector.
\end{enumerate}
Let $W$ be an indecomposable summand of $\cV_s$. By condition $(2)$, one has $W\cong B(\bd,\wt{\lambda})$ with $\bd$ being  nonzero vector with entries consisting of $0,1$. Then it follows from the definition of $ B(\bd,\wt{\lambda})$ that
\[
\dim H^0(W)=\deg W=\sum d_i.
\]
So we have $h^0(V)=n$ and $h^1(V)=0$ for $V\in \cU^+_{L,k+1}$.

\begin{lem}\label{lem: fullrank}
For $V\in\cU^+_{L,k+1}$, the matrix $A(V)$ has rank $k+1$.
\end{lem}
\begin{proof}
From the condition that $\det(V)=L$, it follows that each column of $A(V)$ is a standard basis vector. Moreover the $(i+n)$-th column is the cyclic permutation of the $i$-th column. The full rank property of $A(V)$ follows.
\end{proof}

Recall that $T\subset \ol{T}$ is the subtori defined by $\prod_{i=1}^nt_i =1$ which acts on $C$ as automorphisms. Since $L$ is $T$-equivariant, $\cU^+_{L,k+1}$ inherits a $T$-action.

\begin{lem}\label{lem: binary}
$\cU^+_{L,k+1}$ is an open substack of the stack of vector bundles on $C$ with determinant $L$.
\end{lem}
\begin{proof}
Let $\cV$ be a vector bundle over $C\times S$ where $S$ is a Noetherian base scheme such that $\det(\cV_s)=L$ for $s\in S$. Suppose for $s_0\in S$, $\cV_{s_0}$ satisfies condition $(2)$. We need to show it holds on an open subset containing $s_0$. Condition $(1)$ implies that $\deg(\cV_{s_0}|_{C_i})=1$. Given that, condition $(2)$ holds for $\cV_{s_0}$ if and only if for $i=1,\ldots,n$ there does not exist a nonzero map $\cV_{s_0}\to \OO_{C_i}(-1)$ and there does not exist a nonzero map $\cV_{s_0}\to \OO_C$.  Then the openness of such $s_0$ follows from the semi-continuity theorem.
\end{proof}

The general framework of group action on stacks can be found in \cite{gpstack}. Given an algebraic group $G$ acting on an algebraic stack $\cX$, there is an algebraic quotient stack $\cX/G$. A $G$-orbit (stack) is defined to be the fiber of the canonical map $\cX\to \cX/G$ at a residue gerbe.

\begin{prop}\label{prop: finite orbits}
The stack $\cU^+_{L,k+1}$ has finitely many $T$-orbits.
\end{prop}
\begin{proof}
For simplicity, we set $C=C^n$ and $\wt{C}=C^{rn}$.
Let $V=B(\bd, \wt{\lambda})$ be an indecomposable vector bundle of rank $r$. For $t=(t_1,\ldots,t_n)\in T$,  we get a natural action $\wt{t}: \wt{C}\to \wt{C}$ such that the diagram
\[
\xymatrix{
\wt{C}\ar[r]^{\wt{t}}\ar[d]^{\pi_r} & \wt{C}\ar[d]^{\pi_r} \\
C\ar[r]^t & C
}
\] is Cartesian. Here $\pi_r: \wt{C}\to C$ is a standard $r:1$ cyclic cover.
\[
\wt{t}^*\wt{L}(\bd,\wt{\lambda})=\wt{L}(\bd,\wt{\lambda}^\prime)
\] where
\[
\wt{\lambda}^\prime=\wt{\lambda}\prod_{i=1}^n t_i^{-\sum_{j=0}^{r-1} d_{i+jn}}
\] 
Let $V$ be an arbitrary element in $\cU^+_{L,k+1}$. Then $V=V_1\oplus\ldots\oplus V_p$ with $V_j=B(\bd^j,\wt{\lambda_j})$ indecomposable and $\prod_{j=1}^p\wt{\lambda_j}=\wt{\lambda}$. 
Since for $i=1,\ldots,n$ the determinant condition gives
\[
\sum_{j=1}^p \sum_{\ell=0}^{r_j-1} \bd^j_{i+\ell n}=1,
\]
$\prod_{j=1}^p\wt{\lambda_j}$ is preserved by $\wt{t}^*$. Clearly, the $T$-action preserves the discrete data like rank and $\bd$. For any other $p$-tuple $(\wt{\lambda}^\prime_j)_{j=1}^p$ with $\prod_{j=1}^p\wt{\lambda}^\prime_j=\wt{\lambda}$, the $T$-action is transitive if and only if the system of equations
\[
\wt{\lambda}^\prime_j=\wt{\lambda_j}\prod_{i=1}^n t_i^{-\sum_{\ell=0}^{r_j-1} d_{i+\ell n}}
\] for $j=1,\ldots, p$ has a solution. Since $p\leq k+1\leq n$, this is a underdetermined inhomogeneous linear system (after taking $\log$). It remains to check the matrix $M_{ji}:=\sum_{\ell=0}^{r_j-1} d_{i+\ell n}$ is nondegenerate. By condition $(2)$, $M_{ji}$ is a $p\times n$ binary matrix whose rows are all nonzero. The determinant condition implies that all columns are standard basis vectors, i.e. it has one $1$ and $0$ elsewhere. It follows that $M_{ji}$ has full rank.
\end{proof}

Recall that we have a stack morphism $q$ from $\cN_{L, k}$ to the algebraic stack of vector bundles on $C$ with determinant $L$ that sends the complex $\OO_C^k\to V$ to $V$.
\begin{thm}\label{thm: fibration}
The following properties hold:
\begin{enumerate}
\item[$(1)$] The image of $q$ is contained in the substack $\cU^{+}_{L,k+1}$;
\item[$(2)$] $q$ is $T$-equivariant;
\item[$(3)$] Suppose $v\in  \cU^{+}_{L,k+1}$ is in the image of $q$. Denote by $Tv$ the $T$-orbit of $v$, viewed as a locally closed substack. Then the base change $q^v$
\[
\xymatrix{
\cN^v_{L,k}\ar[r]\ar[d]^{q^v} & \cN_{L,k}\ar[d]^q\\
Tv\ar[r] & \cU^{+}_{L,k+1}
}
\] is smooth and surjective with connected fibers.
\end{enumerate}
\end{thm}
\begin{proof}
Suppose a complex $\OO_C^k\to V$ belongs to $\cN_{L,k}$. By Lemma \ref{lem: binary}, its degree vector $\bd$ must be nonzero and with $0,1$ entries. This proves part $(1)$. Part $(2)$ is obvious since $T$ acts by automorphism of $C$.

Now we prove part $(3)$. The surjectivity of  $q^v$ follows from part $(2)$. 
First we prove that the fiber of $q^v$ is smooth. The proof uses (and only relies on) two  facts: $\cN_{L,k}$ is a $\G_m$-gerbe, and the homotopy fiber product of $q$ and the canonical map from a residue gerbe to $Tv$ is 0-shifted symplectic. We have proved a similar smoothness statement in \cite[Theorem 3.4]{HP-matrix} and the exactly same proof applies to this case as well. We refer the readers to \cite[Section 3]{HP-matrix} for the full details. 

Second we prove the connectedness of the fiber. Let $V$ be the vector bundle representing the point $v\in \cU^+_{L,k+1}$. Let $\phi: \OO_C^k\to V$ be a morphism such that cokernel of $\phi$ is isomorphic to $L$. Injectivity is an open condition. Since the condition that $\cok~ \phi$ is locally free is also an open condition, such morphisms form an open subset $M^s\subset \Hom(\OO_C^k,V)$, therefore is connected. The fiber at $v$ is the quotient stack $\Big[ GL_k\Big\backslash M^s\Big/ \Aut(V)\Big]$. The connectedness of the fiber then follows from the connectedness of $M^s$.

Let  
\[
T_1\to T\to T_2
\]
be the short exact sequence of algebraic groups such that $T_1$ is the stabilizer subgroup of the isomorphism class of $v$.  Then $Tv$ is a $T_1$-gerbe over $T_2$. Recall that $\cN_{L,k}$ is a $\G_m$-gerbe over a smooth scheme $N_{L,k}$ (isomorphic to $G(k,n)$).   Therefore $\cN^v_{L,k}$ a $\G_m$ gerbe over a locally closed subscheme $N^v_{L,k}\subset N_{L,k}$. Given the smoothness and irreducibility of the fiber, it suffices to show the induced scheme morphism  $q: N:=N^v_{L,k}\to T_2$ is flat.   We claim that $N$ is irreducible. Fix $y\in T_2$. Since the fiber is irreducible there an irreducible component $N_1$ of $N$ such that $q^{-1}(y)\subset N_1$. For $t\in T$, by $(2)$ one has
\[
q^{-1}(t\cdot y)=t\cdot q^{-1}(y)\subset N_1.
\] Since $T$ acts transitively on $T_2$, $N=\bigcup_{t\in T}q^{-1}(t\cdot y)=N_1$. Since $T$ acts transitively on $T_2$ and $q$ is equivariant, it follows that all fibers have dimension $\dim(N)-\dim(T_2)$. Hence by \cite[Theorem 3.3.27]{sch-book}, $q$ is flat. So it is smooth.
\end{proof}

\subsection{A geometric construction of the cyclic rank matrix}\label{geom-cyclic-sec}
To finish this section, we give a geometric construction of the cyclic rank matrix (c.f. \cite[Section 3.3]{KLS} or Section \ref{sec: affper}) by identifying $G(k,n)$ with the coarse moduli space of $\cN_{L,k}$.

For each $i\in\Z/n\Z$, set $C_i^\circ$ to be the closure of $C\setminus \{C_{i-1}\cup C_i\}$.  For $i\leq j\leq i+n-2$, set
\[
C^\circ_{ij}=C_i^\circ \cap C^\circ_{i+1}\cap \ldots\cap C^\circ_j.
\]
Note that for $j-i=n-2$ 
\[
C^\circ_{ij}=C_{i-2}\cap C_{i-1}. 
\] When $j-i<n-2$, 
\[
C^\circ_{ij}=C_{j+1}\cup C_{j+2}\cup\ldots\cup C_{i-2}. 
\] 
Given a vector bundle $V$ on $C$, we denote 
\[
H^0_{[i-1,j]}(V):=\ker\Big(H^0(V)\to H^0(V|_{C^\circ_{ij}}) \Big).
\] We set
\[
h_{ij}(V)=\dim\Big(H^0_{[i-1,j]}(V)\Big).
\] 
Note that since $C^\circ_{i,i+n-2}$ is a point, we have $\dim H^0(V|_{C^\circ_{i,i+n-2}})=\rk V$. Hence,
\begin{equation}\label{hij-point-inequality}
h_{i,i+n-2}\ge h^0(V)-\rk V.
\end{equation}

Suppose that $V$ is an extension
\[
\OO_C\ot W\to V\to L
\] such that $\dim(W)=k$ and $\deg(L|_{C_i})=1$ for $i\in\Z/n\Z$. We choose local coorindate $x_i$ on $C_i$ such that $C_{i-1}\cap C_i=\{x_{i-1}=\infty\}=\{x_i=-1\}$. Without loss of generality, we assume that $L\cong \OO_C(D)$ for 
\[
D=\bigcup_{i\in\Z/n\Z}\{x_{i}=0\}.
\]
The ring of regular functions on $C\setminus D$ is the subalgebra of $\prod_{i\in \Z/n\Z}\C[x_i^{-1}]$ generated by
\[
s_{i}:=(1+x_{i-1}^{-1})e_{i-1}-x_i^{-1}e_{i}
\] where $e_i$ is the $i$-th idempotent. In Section \ref{sec: standard}, we have shown that $\{s_i|i\in \Z/n\Z\}$ gives a basis of $H^0(L)$. Clearly, $s_i$ is supported on $C_{i-1}\cup C_i$. The connecting morphism of the short exact sequence defining $V$ is a linear map
\[
M_V: H^0(L)\to W.
\] When $V$ is fixed we write $M$ for $M_V$ for simplicity.
Following \cite{KLS},  we define a  matrix $r(M):=\Big(r_{ij}(M)\Big)_{i,j\in\Z}$ by
\begin{equation}\label{r(M)}
r_{ij}(M)=
\begin{cases}
\dim\Big(\sspan\Big(M(s_i),\ldots,M(s_j)\Big)\Big) & i\leq j\\
j-i+1 & j<i
\end{cases}
\end{equation} where
$s_i=s_{i\mod n}$. 
 By \cite[Lemma 5.2]{KLS}, $r_{ij}(M)$ satisfies the axioms of the \emph{cyclic rank matrix} (see \cite[Corollary 3.12]{KLS}, c.f. Proposition \ref{prop: cyclic}).

The following Lemma gives a geometric interpretation of cyclic rank.
\begin{lem}\label{lem: r(M)}
Let 
\[
\OO_C\ot W\to V\to L
\] be an extension such that $M: H^0(L)\to W$ is surjective.
Then we have an equality
\[
r_{ij}(M)=j-i+1-h_{ij}(V)
\] for $i\leq j\leq i+n-2$.
\end{lem}
\begin{proof}
There is a morphism of long exact sequences
\[
\xymatrix{
W\ar[r]\ar[d]^= & H^0(V)\ar[d]\ar[r] & H^0(L) \ar[r]^M\ar[d]^\iota &W\ar[d]\\
W\ar[r] & H^0(V|_{C^\circ_{ij}})\ar[r]& H^0(L|_{C^\circ_{ij}})\ar[r] & H^1(C^\circ_{ij},\OO_C\ot W)=0.
}
\] 
It follows that 
\[
H^0_{[i-1,j]}(V)\cong \ker\Big(H^0(V)/W\to H^0(V|_{C^\prime_{ij}})/W\Big).
\]
The natural map 
\[
H^0_{[i-1,j]}(V)\to \ker(\iota)
\] is injective and factors through $\ker(\iota)\cap \ker(M)$.  It is surjective by the vanishing of $H^1(C^\circ_{ij},\OO_C\ot W)$ since $C^\circ_{ij}$ is a chain of rational curves.
Applying $\iota$ on the basis $s_1,\ldots,s_n$, we get that 
\[
\dim\Big(\ker(\iota)\Big)=j-i+1.
\] Then the equality follows from the exactness.
\end{proof}

 For a line bundle $\wt{L}$ on $\wt{C}:=C^{rn}$ with the degree vector 
\[
\bd=(d_1,\ldots,d_{rn}).
\] set $V=(\pi_r)_*\wt{L}$. For $i\leq j\leq i+n-1$, set 
\begin{equation}\label{sij-eq}
s_{ij}(\bd):=\max(d_{i-1}+d_i+\ldots +d_j-1,0),
\end{equation}
where the indices are modulo $rn$. Below we will use the same function $s_{ij}$ on infinite sequences of $0$'s and $1$'s indexed by $\Z$.

\begin{lem}\label{lem: formula of h(M)}
In the above situation we have 
\[
h_{ij}(V)=\sum_{\ell=0}^{r-1} s_{i+n\ell,j+n\ell}(\bd), 
\] for $i\leq j\leq i+n-2$.
\end{lem}
\begin{proof}

Note that $H^0(C,V)=H^0(\wt{C}, \wt{L})$. Denote by $\wt{C}^\circ_{ij}$ the preimage of $C^\circ_{ij}$.
We then have
\[
\overline{\wt{C}\setminus \wt{C}^\circ_{ij}}=\bigsqcup_{k=0}^{r-1} \wt{C}_{[i-1+nk,j+nk]}
\]
where $\wt{C}=\bigcup_{k=1}^{rn} C_k$ and $\wt{C}_{[i,j]}=\bigcup_{k=i}^{j} C_k$. By base change we can identify $H^0_{[i-1,j]}(V)$ with $\bigoplus_{k=0}^{r-1} W_{[i-1+nk,j+nk]}$ where $W_{[i-1+nk,j+nk]}\subset H^0(\wt{C}_{[i-1+nk,j+nk]},\wt{L})$ consists of sections that vanish at $x_{i-1+nk}=\infty$ and $x_{j+nk}=-1$.
It is easy to check that 
\[
\dim W_{[i-1+nk,j+nk]}= s_{i+nk,j+nk}(\bd).
\]
\end{proof}

Now let $V$ be an arbitrary vector bundle, and let $A(V)$ be the corresponding matrix.

\begin{lem}\label{hijV-lem}
For $i\leq j\leq i+n-2$, one has 
$$h_{ij}(V)=\sum_{\ell=1}^r s_{ij}(w_\ell)$$
where $w_1,\ldots,w_r$ are the rows of the matrix $A(V)$ (so $r$ is the rank of $V$). 
\end{lem}

\begin{proof}
For an indecomposable $V$, this follows from Lemma \ref{lem: r(M)} and the fact that
$w_{i+1}=\tau w_i$ in this case. The case of an arbitrary $V$ follows due to the way the matrix $A(V)$ is formed.
\end{proof}

We extend the range of $h_{ij}(V)$ for any vector bundle $V$ of rank $k+1$ by setting
\begin{equation}\label{h(V)}
h_{ij}(V)=h_{ij}(A):=
\begin{cases}
\sum_{\ell=1}^{k+1} s_{ij}(w_\ell), & i\leq j\leq i+n-2,\\
0, & j<i,\\
j-i+1-k, & j\geq i+n-1,
\end{cases}
\end{equation} 
where $A=A(V)$ and $w_1,\ldots,w_{k+1}$ are the rows of $A$. By Lemma \ref{hijV-lem}, this agrees with the previous definition
in the range $i\le j\le i+n-2$.


Now we define the matrix $r(V)$ by
\[
r_{ij}(V)=j-i+1-h_{ij}(V).
\]

\begin{lem}\label{rM-rV-lem}
For an extension
\[
\OO_C\ot W\to V\to L
\] that induces a surjective map $M: H^0(L)\to W$, one has
\begin{equation}\label{rM-rV-eq}
r(M)=r(V).
\end{equation}
\end{lem}

\begin{proof}
This follows from  Lemma \ref{lem: r(M)}.
\end{proof}

\begin{rmk}
One can show that for a vector bundle $V$ of rank $k+1$ with $h^0(V)=n$, the formula $h_{ij}(V)=\sum_{\ell=1}^{k+1} s_{ij}(w_\ell)$ actually holds
for all $i\le j$, i.e., this expression is equal to $j-i+1-k$ for $j\ge i+n-1$. Indeed, the inequality \eqref{hij-point-inequality} implies that $h_{i,i+n-2}\ge n-1-k$.
This implies that the rows matrix $A$ cannot have $n$ consecutive zeros, and this leads to the claimed assertion.
\end{rmk}

The assignment $V\mapsto h(V)$ defines a matrix valued function on $\cU^+_{L,k+1}$. For matrices $A$ and $B$ we define a partial order $A\geq B$ if and only if for all $i,j$ $A_{ij}\geq B_{ij}$.
\begin{prop}\label{prop: h semicont}
The matrix valued function $h(V)$ is upper semi-continuous on $\cU^{+}_{L,k+1}$.
\end{prop}
\begin{proof}
It suffices to check that $h_{ij}(V)$ is upper semi-continuous  for $0\leq j-i\leq n-2$. Let $S\to \cU^+_{L,k+1}$ be a morphism from an affine scheme $S$. Let $\cV_S$ be the pull back of the universal family over $C\times S$. We push forward the sheaf morphism $\cV_S\to \cV_S\ot p_1^*\OO_{C^\circ_{ij}}$ along the projection $p_2: C\times S\to S$ to obtain a sheaf morphism 
\[
\phi: (p_2)_*\cV_S\to (p_2)_*\Big(\cV_S\ot p_1^*\OO_{C^\circ_{ij}}\Big).
\] 
We claim this is a morphism of vector bundles. For any $s\in S$, one has $\dim H^0(\cV_s)=n$ by condition $(2)$. Since $C^\circ_{ij}$ is a chain (not cycle) of rational curves $\cV_s|_{C^\circ_{ij}}$ splits. Moreover the splitting type of $\cV_s|_{C_i}$ if independent of $s$, therefore $\dim H^0(\cV_s|_{C^\circ_{ij}})$ is independent of $s$. The claim then follows from the formal function theorem.
Then
\[
h_{ij}(\cV_s)=n-\rk(\phi_s),
\] is upper semi-continuous since $\rk(\phi_s)$ is lower semi-continuous.
\end{proof}

\section{Affine permutation and cyclic rank matrix}\label{sec: affper}
\subsection{Affine permutations}
Given $k<n\in \N$, let $B(k,n)$ be the set of \emph{bounded affine permutations of average $k$} which are bijections $f:\Z\to \Z$ such that
\begin{enumerate}
\item[$(1)$] $f(i+n)=f(i)+n$;
\item[$(2)$] $i\leq f(i)\leq n+i$;
\item[$(3)$] $\sum_{i=1}^n (f(i)-i)=kn$.
\end{enumerate}
A bijection satisfying $(1)$ is called an \emph{affine permutation}. The group of all affine permutations (of period $n)$ is denoted by $\wt{S}_n$. An affine permutation satisfying condition $(2)$ is called \emph{bounded}. The positive integer $k$ in condition $(3)$ is called the \emph{ball number} or the \emph{average} of an affine permutation. The subset of $\wt{S}_n$ satisfying $(3)$ is denoted by $\wt{S}^k_n$.  We denote by $\wt{S}^{k,+}_n\subset \wt{S}^k_n$ the subset satisfying the condition $f(i)\geq i$ and the subset  $\wt{S}^{k,++}_n\subset \wt{S}^{k,+}_n$ those satisfying $f(i)>i$. Note that since condition $(3)$ gives a upper bound $f(i)-i\leq kn$,   $\wt{S}^{k,+}_n$ is a finite set.

\begin{lem}\label{lem:orbit}
Given $f\in\wt{S}^{k,++}_n$, $\Z_f:=\{f^d|d\in\Z\}\subset\wt{S}_n$ is an infinite cyclic group.
The action of $\Z_f$ on $\Z$ has exactly $k$ orbits.
\end{lem}
\begin{proof}
Condition $(3)$ implies that 
\[
\lim_{m\to +\infty} \frac{1}{2mn} \sum_{i=-mn+1}^{mn} (f(i)-i)= k.
\]
Fix an orbit $\alpha$ of $\Z_f$. We have
\[
f((m-1)n)+(m-1)n-1\leq \sum_{i\in \alpha\cap [-mn+1,mn]} (f(i)-i)\leq f(mn)+mn-1. 
\] Applying the boundness  condition, we obtain an equality 
\[
2(m-1)n-1\leq \sum_{i\in \alpha\cap [-mn+1,mn]} (f(i)-i) \leq 2mn+kn-1.
\]
Taking the limit we obtain
\[
\lim_{m\to +\infty} \frac{1}{2mn}\sum_{i\in \alpha\cap [-mn+1,mn]} (f(i)-i)=1.
\]
Combining with condition $(3)$ we prove the desired statement.
\end{proof}

Let $\alpha\subset \Z$ be an orbit of $\Z_f$. The \emph{characteristic vector} $w^\alpha=(w^\alpha_i)_{i\in\Z}$ of $\alpha$ is an infinite vector such that
$w^\alpha_i=1$ if $i\in \alpha$ and zero otherwise. Given an affine permutation $f\in\wt{S}^{k,++}_n$, the \emph{characteristic matrix} $A_f$ of $f$ is the $k\times \infty$-matrix $A$ whose 
rows are characteristic vectors of all the orbits of $\Z_f$.  It is defined uniquely up to row permutations.

Denote by $\ol{f} \in S_n$ the underlying finite permutation of $f$ modulo $n$ and by $\ol{\alpha}$ the finite cycle underlying $\alpha$. Let $|\ol{\alpha}|$ be the length of the cycle. We represent the orbit $\alpha$ by an increasing monotone sequence
\[
\alpha=\Big( \ldots \alpha_{i-1}<\alpha_i<\alpha_{i+1}<\ldots\Big).
\]

\begin{lem}\label{lem: length}
For any $i\in\Z$, $n\Big| \alpha_{i+|\ol{\alpha}|}-\alpha_i$ and the ratio, denoted by $\omega(\alpha)$ is independent of $i$. We call it the \emph{period} of $\alpha$.
\end{lem}
\begin{proof}
The first claim holds obviously. Observe that $\alpha_{i+|\ol{\alpha}|}=f^{|\ol{\alpha}|}(\alpha_i)$. Suppose that 
\[
f^{|\ol{\alpha}|}(\alpha_i)=\alpha_i+\omega n.
\] By applying $f$ on both side and apply condition $(1)$, we get
\[
f^{|\ol{\alpha}|}(\alpha_{i+1})=f^{|\ol{\alpha}|+1}(\alpha_i)=f(\alpha_i)+\omega n=\alpha_{i+1}+\omega n.
\] The second claim then follows from induction.
\end{proof}
 
We denote by $s_+$ the permutation $s_+(i)=i+1$.
\begin{lem}\label{lem: period}
Let $\alpha$ be an orbit of $\Z_f$. Then $s_+^n(\alpha)$ is also an orbit. Moreover,
\[
\omega(\alpha)=\min\{d>0\big| s_+^{dn}(\alpha)=\alpha\}.
\]
\end{lem}
\begin{proof}
This is clear.
\end{proof}
If we denote by $s_i$ the transposition $(i, i+1)$, then the subgroup $\wt{S}^0_n\subset \wt{S}_n$ is a Coexter group generated by $s_0,s_1,\ldots,s_{n-1}$, therefore admits a Bruhat order. Since $\wt{S}^k_n$ can be identified with the coset $s_+^k \wt{S}^0_n$, it becomes a poset under the Bruhat order of $\wt{S}^0_n$. To be more precise, for $g,f\in \wt{S}^k_n$ we write $g\precdot f$ if $g$ is obtained from $f$ by swapping the values of $i+kn$ and $j+kn$ for all $k$, where $i<j$ and $f(i)>f(j)$. The Bruhat order is its transitive closure. 
The element $s_+^k$ is minimal under the Bruhat order.

\subsection{Cyclic rank matrix and its generalization}\label{cyclic-rank-sec}
We  recall the definition of a cyclic rank matrix and its generalization following \cite{KLS}. 
\begin{prop}\cite[Corollary 3.10, 3.12]{KLS}\label{prop: cyclic}
Fix $0<k<n$. 
The set of $\infty\times \infty$ integer matrices $r_{ij}$ satisfying condition $C_1,\ldots,C_4$ is in one to one correspondence with $\wt{S}^k_n$ under the assignment $r_{ij}\mapsto f$ where $f(i)=j$ if and only if 
\[
r_{ij}=r_{i+1,j}=r_{i,j-1}>r_{i+1,j-1}.
\]
\begin{enumerate}
\item[$(C_1)$] for each $i$, there exists $m_i$ such that $r_{ij}=j-i+1$ for all $j\leq m_i$;
\item[$(C'_1)$] $r_{ij}=j-i+1$ when $j<i$;
\item[$(C_2)$] for each $i$, there exists $n_i$ such that $r_{ij}=k$ for all $j\geq n_i$;
\item[$(C'_2)$] $r_{ij}=k$ when $j-i\geq n-1$;
\item[$(C_3)$] $r_{ij}-r_{i+1,j}, r_{ij}-r_{i,j-1}\in \{0,1\}$;
\item[$(C_4)$] if $r_{i+1,j-1}=r_{i+1,j}=r_{i,j-1}$ then $r_{i+1,j-1}=r_{ij}$;
\item[$(C_5)$] $r_{i+n,j+n}=r_{ij}$.
\end{enumerate}
We call such a matrix a \emph{generalized cyclic rank matrix}.
Moreover, under this correspondence the subset consisting of matrices satisfying the additional conditions $C'_1,C'_2,C_5$ is bijective to $B(k,n)\subset \wt{S}^k_n$. We call such a matrix a \emph{cyclic rank matrix}.
\end{prop} 
Equivalently, we may state the above conditions for the matrix $h_{ij}=j-i+1-r_{ij}$ and get:
\begin{enumerate}
\item[$(H_1)$] for each $i$, there exists $m_i$ such that $h_{ij}=0$ for all $j\leq m_i$;
\item[$(H'_1)$] $h_{ij}=0$ when $j<i$;
\item[$(H_2)$] for each $i$, there exists $n_i$ such that $h_{ij}=j-i+1-k$ for all $j\geq n_i$;
\item[$(H'_2)$] $h_{ij}=j-i+1-k$ when $j-i\geq n-1$;
\item[$(H_3)$] $h_{ij}-h_{i+1,j}, h_{ij}-h_{i,j-1}\in \{0,1\}$;
\item[$(H_4)$] if $h_{i+1,j-1}+1=h_{i+1,j}=h_{i,j-1}$ then $h_{i+1,j-1}+2=h_{ij}$;
\item[$(H_5)$] $h_{i+n,j+n}=h_{ij}$.
\end{enumerate}

We write $r(f)$ and $h(f)$ the $r$-matrix and $h$-matrix associated to the affine permutation $f\in \wt{S}^k_n$. 
For two matrices $r, s$ satisfying conditions $C_1, C_2, C_3,C_4$, we write $r\geq s$ if for any $i,j\in\Z$ $r_{ij}\geq s_{ij}$. This order is related to the Bruhat order by the following proposition.
\begin{prop}\cite[Proposition 3.11]{KLS}\label{prop: order on r} 
Let $f,g\in \wt{S}^k_n$ and let $r$ and $s$ be the corresponding matrices. Then  $f\preceq g$ under the Bruhat order if and only if $r\geq s$.
\end{prop}

For $0<r<N$ let $M=[v_1,\ldots,v_N]$  be an $r\times N$ complex matrix with column vectors $v_1,\ldots,v_N$. Define the \emph{cyclic rank} function $f_M:\Z\to \Z$  by
\[
f_M(i)=\min\{j\geq i\big| v_i\in  \sspan(v_{i+1},\ldots,v_j)\}
\] where the column index is considered modulo $N$.

\begin{prop}\cite[Proposition 3.15, Lemma 4.2]{KLS}\label{prop: cyclic rk}
If $M$ has rank $r$ then $f_M\in B(r,N)$. Furthermore, $f_M$ corresponds to the cyclic rank matrix $r(M)$ under the bijection of Proposition \ref{prop: cyclic}.
\end{prop}

\section{Torus orbits of moduli stacks of vector bundles and extensions} \label{sec: main}
We prove the two major technical results of the paper in this section. The first one is a combinatorial description of the $T$-orbits of $\cU^+_{L,k+1}$ and the second one is a classification theorem of extensions of the form $\OO_C^k\to V\to L$.

\subsection{Counting torus orbits of $\cU^+_{L,k+1}$}
Let $A=[v_1,\ldots,v_N]$ be a $(k+1)\times N$ matrix whose columns $v_i$ are standard basis vectors and all rows are nonzero (a standard basis vector is a vector with exactly one entry $1$ and $0$'s elsewhere). By an abuse of notation we denote by the same symbol $A$ the $(k+1)\times \infty$-matrix that is  the $N$-periodic extension of $A$ in columns.  We define a function $f_A:\Z\to \Z$ by 
\[
f_A(i)=\min\{j\ge i\big|, v_j=v_{i-1}\}.
\]  Clearly, $f_A$ is a bijection and is invariant under row permutations of $A$.
\begin{thm}\label{thm: Vtof} 
The assignment $\phi: V\mapsto f_A$ for $A=A(V)$ defined in Section \ref{sec: vb} satisfies 
\begin{enumerate}
\item[$(1)$]
$\phi$ induces a map $[\phi]: \cU^+_{L,k+1}/T\to \wt{S}^{k,+}_n$;
\item[$(2)$]  $[\phi]$ is a bijection.
\end{enumerate}
\end{thm}
\begin{proof} We prove part $(1)$ first. Recall that $A(V)$ has size $(k+1)\times hn$ and we extend it $hn$-periodically to a $(k+1)\times\infty$-matrix.
It is obvious that $f_A(i)\geq i$. For $V=V_1\oplus\ldots \oplus V_p$ where $V_i$ are indecomposables of rank $r_i$, the cyclic column permutation
\[
(v_1,\ldots,v_{hn})\mapsto (v_{n+1},\ldots,v_{hn},v_1,\ldots,v_n) 
\] where $h=\lcm(r_1,\ldots,r_p)$
differs from $A$ by a row permutation. This leads to the property $f_A(i+n)=f_A(i)+n$. 
By Lemma \ref{lem: fullrank} and Proposition \ref{prop: cyclic rk}, $f_A\circ s_+\in B(k+1,hn)$, so we have 
\[
\sum_{i=1}^{hn} (f_A(i)-i+1)=(k+1)hn.
\] Hence, $\sum_{i=1}^{hn} (f_A(i)-i)=khn$. Then by (column) $n$-periodicity, we get that 
$$\sum_{i=1}^{hn}(f_A(i)-i)=h\cdot \sum_{i=1}^n (f_A(i)-i)$$ and deduce 
that $\sum_{i=1}^{n} (f_A(i)-i)=kn$. Since $A(V)$ only depends on the degree vector of the indecomposable summands of $V$, it is preserved by the $T$-action. So the map $V\mapsto f_{A(V)}$ is well defined on the $T$-orbits of $\cU^+_{L,k+1}$.

To prove part $(2)$,  we are going to construct an inverse map $\wt{S}^{k,+}_n\to \cU^{+}_{L,k+1}/T$. Let $A=A(V)$ for some $V\in \cU^+_{L,k+1}$ and $f=f_{A}$. Recall that we can associate to $f$ the characteristic matrix $A_{fs_+}$. It is clear that for such $A$
\[
A_{f_As_+}=A
\] up to row permutations. It remains to check that for any $f\in \wt{S}^{k,+}_n$, $A_{fs_+}=A(V)$ for some vector bundle $V$.

Denote by $\mathfrak{C}_{fs_+}$ the set of $\Z_{fs_+}$ orbits. By Lemma \ref{lem:orbit} it has cardinality $k+1$. By Lemma \ref{lem: period}, the infinite cyclic group $\{s_+^{dn}\big|d\in\Z\}$ acts on $\mathfrak{C}_{fs_+}$. We call  $\alpha, \beta$   equivalent if $s_+^{nd}(\alpha)=\beta$ for some $d$. Denote by $p(f)$ the number of orbits of this group action. The equivalence class of $\alpha$ has cardinality $\omega(\alpha)$. We associate to each equivalence class $[\alpha]$ a (numerical type of) indecomposable vector bundle $[V_\alpha]$ with rank $\omega(\alpha)$ and degree $|\ol{\alpha}|$. The degree vector of $[V_\alpha]$ is set to be $w^{\alpha}$ (truncate its fundamental domain under translation). Recall that $w^{\alpha}$ and $w^{s_+^n(\alpha)}$ determine the identical numerical type of indecomposable bundles. Finally, we have $A_{fs_+}=A(V)$ for 
\[
V=\bigoplus_{[\alpha], \alpha\in \mathfrak{C}_f} V_\alpha
\] by choosing $V_\alpha$ with the prescribed numerical type $[V_\alpha]$.  Since the $T$-action preserves the numerical type of $V$ and the action is transitive (Proposition \ref{prop: finite orbits}) we get a well defined inverse map.
\end{proof}

\begin{rmk}
A by-product of Theorem \ref{thm: Vtof} is a nice correspondence between combinatorial invariants of affine permutations and algebraic geometric invariants of vector bundles on Kodaira cycles. Given $f\in \wt{S}^{k,+}_n$, the corresponding vector bundle $V_f$ has rank $k+1$ and degree $n$. The number of indecomposable summands $p$ of $V$ is equal to $p(f)$, i.e. the number of orbits in $\mathfrak{C}_{fs_+}$ under the action of the cyclic group $\{s_+^{nd}|d\in\Z\}$. And indecomposable summands $V_\alpha$ are in one to one correspondence with the orbits $[\alpha]$ mentioned above. The rank of $V_\alpha$ is equal to the period $\omega(\alpha)$ and degree of $V_\alpha$ is equal to $|\ol{\alpha}|$. The degree vectors of $V_\alpha$ are encoded in the rows of  the matrix $A(V)$. Later we will show that the length $\ell(f)$ is equal to the dimension of $\End(V)$ minus $p(f)$(see Proposition \ref{prop: end}). 
\end{rmk}

\begin{prop}\label{cor: U++}
We set 
\[
\cU^{++}_{L,k+1}:=\bigsqcup_{f\in B(k,n)} T_f, 
\] where $T_f:=[\phi]^{-1}(f)$ is the $T$-orbit that corresponds to $f$. Then 
$\cU^{++}_{L,k+1}$ is an open substack of $\cU^{+}_{L,k+1}$. 
\end{prop}
\begin{proof}
Let $v, w$ be two $T$-orbits in $\cU^+_{L, k+1}$. We say $v\prec w$ if $v$ is contained in the closure of $w$. Suppose $\phi(V)=f\in \wt{S}^{k,+}_n$. We denote by $[f, +\infty), (-\infty,f]$ the inteverals in $\wt{S}^{k,+}_n$ under the Bruhat order.  We claim that $\phi: \cU^{+}_{L,k+1}\to \Big(\wt{S}^{k,+}_n, \preceq\Big)$ is upper semi-continuous where $\wt{S}^{k,+}_n$ is equipped with the order topology. First note that the Bruhat order is opposite to the order of $r$-matrix by Proposition \ref{prop: order on r}. Since $r_{ij}=j-i+1-h_{ij}$, the claim then follows from Proposition \ref{prop: h semicont}. Since $\wt{S}^{k,+}_n$ has finitely many $T$-orbits by Proposition \ref{prop: finite orbits}, the upper semi-continuity implies that $\phi^{-1}(-\infty, f]$ is open and $\phi^{-1}[f,+\infty)$ is closed.  The subset $B(k,n)\subset \wt{S}^k_n$ is a lower order ideal, i.e. if $g\preceq f$ for $f\in B(k,n)$ then $g\in B(k,n)$ (\cite[Lemma 3.6]{KLS}). Then 
\[
\cU^{++}_{L,k+1}=\bigcup_{\text{$f$ maximal in $B(k,n)$}} \phi^{-1}(-\infty, f]
\] is open.
\end{proof}

\subsection{Classification of extensions}

\begin{thm}\label{thm: extension}
The stack morphism $q: \cN_{L,k}\to \cU^+_{L,k+1}$ factors through the substack $\cU^{++}_{L,k+1}$.  Moreover 
\[
[\phi]: \cU^{++}_{L, k+1}/T\to B(k,n) 
\]
defined in Theorem \ref{thm: Vtof}, is  an anti-isomorphism of posets, where the partial order on $\cU^+_{L,k+1}/T$ (and restricted to $\cU^{++}_{L,k+1}/T$) is given by the adjacency order of the $T$-orbits.
\end{thm}

\begin{proof}
Let $\OO_C\ot W\to V\to L$ be an extension in $\cN_{L,k}$ for a vector space $W$ of dimension $k$, and let $M: H^0(L)\to W$ be the corresponding  surjection.
We know that $r(V)=r(M)$ (see \eqref{rM-rV-eq}), where $r(M)$ is the cyclic rank matrix of type $(k,n)$ associated with $M$ (see \eqref{r(M)}). Let $f_M$ be the corresponding element of $B(k,n)$ via the bijection
of Proposition \ref{prop: cyclic}.
On the other hand, we have the matrix $A=A(V)=(v_1,\ldots,v_{hn})$ and the corresponding affine permutation
$f=f_A\in \wt{S}^{k,+}_n$ where
\[
f_A(i)=\min\{j\ge i\big|, v_j=v_{i-1}\}.
\] 
We claim that $f=f_M$, in particular, $f\in B(k,n)$.


Recall (see Sec.\ \ref{geom-cyclic-sec}) that $r=r(M)=r(V)$ satisfies $r_{ij}=j-i+1-h_{ij}$, where $h_{ij}$
is given by
\begin{enumerate}
\item[$(1)$] $h_{ij}=0$ for $j<i$;
\item[$(2)$] $h_{ij}=j-i+1-k$ for $j-i\geq n-1$;
\item[$(3)$] for $i\leq j< i+n-1$, $h_{ij}=\sum_{\ell=1}^{k+1} s_{ij}(w_\ell)$;
\item[$(4)$] $h_{i+n,j+n}=h_{ij}$.
\end{enumerate}
where $w_1,\ldots,w_{k+1}$ are the (infinite $nh$-periodic) rows of $A$, and
$s_{ij}(\cdot)$ is defined by \eqref{sij-eq}.

By definition, if $f(i)=j$ then $v_j=v_{i-1}=e_\ell$ for some $\ell$ and the truncated row
$w_{\ell}^{[i,j-1]}$ is zero. This implies that $h_{ij}=h_{i,j-1}+1=h_{i+1,j}+1=h_{i+1,j-1}+1$. In particular, $h_{i+1,j-1}>h_{ij}-2$.
By definition of the bijection of Proposition \ref{prop: cyclic} this implies that $f_M(i)=j$. Hence, $f=f_M$. 


Conversely, starting with $f\in B(k,n)$, by the result of \cite[Sec.\ 5.2]{KLS}, we can find a surjection $M:H^0(L)\to W$ such that 
$r(M)$ corresponds to $f$ via the bijection of Proposition \ref{prop: cyclic}.
By Serre duality, $M$ determines an element in $\Ext^1(L,\OO_C\ot W)$, and in particular a point $\OO_C\ot W\to V\to L$ in $\cN_{L,k}$. 
As we proved above $f=f_{A(V)}$, so this proves the surjectivity of the map $\cN_{L,k}\to \cU^{++}_{L,k+1}$.

To prove the second part, it suffices to show that
\[
\ol{T_f}=\bigsqcup_{f\preceq g \in B(k,n)} T_g,
\] where the closure is taken in $\cU^{++}_{L,k+1}$.
The left hand-side is contained in the right hand-side since $\bigsqcup_{f\preceq g} T_g=\phi^{-1}[f,+\infty]$ is closed. Recall that in part $(3)$ of Theorem \ref{thm: fibration}, we prove that $q^v: \cN^v_{L,k+1}\to Tv$ is surjective. For simplicity we set $q^f=q^v: \cN^f\to T_f$ for any $v\in\phi^{-1}(f)$. 
Given $f\prec g\in B(k,n)$,  to check $T_g\subset \ol{T_f}$ it suffices to show that $\cN^g\subset \ol{\cN^f}$.
Since this is a topological property, we may instead verify $X_g\subset \ol{X_f}$ where $X_g,X_f$ are the coarse moduli spaces of $\cN^g, \cN^f$. By sending extensions $\OO_C\ot W\to V\to L$ to the linear map  $M: H^0(L)\to W$ we identify $X_f$ with the open positroid variety with affine permutation $f$ (see Definition \ref{def: positroid}). The property $X_g\subset \ol{X_f}$ for $f\prec g$ has been proved for positroid varieties (c.f. \cite[Section 5.2]{KLS}).
\end{proof}

In the next theorem, we are going to relate the length of a bounded affine permutation to the dimension of the endomorphism ring of the corresponding vector bundle. This result has independent interest  in Poisson geometry since it determines the dimension of $T$-leaves. In the algebraic geometry counterpart, $\ell(f)$ measures the failure of reductiveness of the automorphism group of the vector bundle $V_f$.

\begin{lem}\label{lem: H0}
Let $V=B(\bd,1,\wt{\lambda})$ be an indecomposable bundle of rank $r$ on $C=C^n$. A nonnegative subsequence $\bd_I:=\{d_i|i\in I=[i^--1, i^+]\}\subset \bd$ is called maximal if $d_{i^--2}<0$ and $d_{i^++1}<0$. Here the index $i$ is taken in $\Z/rn\Z$. Let $\{I_1,\ldots,I_a\}$ be the set of all maximal nonnegative sequences of $\bd$. We define a function
\[
\theta(\bd,\wt{\lambda})=\begin{cases}
1 & \text{if $\bd={\bf{0}}$ and $\wt{\lambda}=1$},\\
0 & \text{if $\bd={\bf{0}}$ and $\wt{\lambda}\neq 1$},\\
\sum_{\ell=1}^a s_{i^-_\ell i^+_\ell}(\bd_{I_\ell}) & \text{otherwise},
\end{cases}
\]
where we use the function $s_{ij}(\cdot)$ defined by \eqref{sij-eq}.
Then
\[\dim\Big(H^0(C, B(\bd,1,\wt{\lambda}))\Big)=\theta(\bd,\wt{\lambda}).
\]
\end{lem}
\begin{proof}
This is a special case of \cite[Lemma 3.23]{BBDG}.
\end{proof}

\begin{lem}
For $i=1,2$, let $V_i=B(\bd^i,1, \wt{\lambda_i})$ be two indecomposable vector bundles of rank $r_i$ on $C=C^n$. Set $h=\lcm(r_1,r_2)$, $g=\gcd(r_1,r_2)$ and define $\wt{\bd^i}\in \Z^{hn}$, for $i=1,2$, by the formula
\[
\wt{\bd^i}=(\bd^i,\bd^i,\ldots,\bd^i)
\] with $h/r_i$ repeated copies of $\bd^i$. Recall that for $\wt{\bd}\in \Z^{hn}$, we set $\tau(\wt{\bd})_i:=\wt{\bd}_{i+n}$. Then
\[
\dim\Big(\Hom(V_1,V_2)\Big)=\sum_{\ell=0}^{g-1} \theta\Big(\tau^{\ell}(\wt{\bd_2})-\wt{\bd_1},\frac{\wt{\lambda_2}^{r_2/g}}{\wt{\lambda_1}^{r_1/g}}\Big).
\]
\end{lem}
\begin{proof}
Let $\pi_i: C^{r_in}\to C^n$ be the standard cyclic cover for $i=1,2$. Then we have line bundles $\wt{L_i}=\wt{L}(\bd^i,\wt{\lambda_i})$ on $C^{r_in}$ such that $V_i=(\pi_i)_*\wt{L_i}$. Consider the Catesian diagram
\[
\xymatrix{
\wt{C}\ar[r]^{p_1}\ar[d]^{p_2} & C^{r_1n}\ar[d]^{\pi_1}\\
C^{r_2n}\ar[r]^{\pi_2} & C^n
}
\] and set $\wt{\pi}=\pi_2\circ p_2=\pi_1\circ p_1$.
By \cite[Lemma 3.24]{BBDG}, 
\[
V_1^\vee\ot V_2\cong \wt{\pi}_*\Big(p_1^*\wt{L_1}^{\vee}\ot p_2^*\wt{L_2}\Big).
\] We then have
\[
\dim\Big(\Hom(V_1,V_2)\Big)=\dim H^0(\wt{C}, p_1^*\wt{L_1}^{\vee}\ot p_2^*\wt{L_2}).
\]
By \cite[Proposition 3.25]{BBDG}, 
\[
\wt{C}=\bigsqcup_{\ell=0}^{g-1} C^{hn}_{(\ell)},
\] where $C^{hn}_{(\ell)}$ is isomorphic to $C^{hn}$ for $\ell=0,\ldots,g-1$. Moreover $p_1|_{C^{hn}_{(\ell)}}$ is the standard cyclic cover $\pi_{h/r_1}: C^{hn}\to C^{r_1n}$, and $p_2|_{C^{hn}_{(0)}}$ is the standard cyclic cover $\pi_{h/r_2}: C^{hn}\to C^{r_2n}$. For $\ell>0$, 
$p_1|_{C^{hn}_{(\ell)}}$ is the standard cyclic cover precomposed by $\tau^{\ell }$ where $\tau: C^{hn}\to C^{hn}$ is the automorphism that sends the $j$-th component to the $(j+n)$-th component. For simplicty, we set $p^\ell_i=p_i|_{C^{hn}_{(\ell)}}$. We have 
\[
(p_1^\ell)^*\wt{L_1}^\vee\ot (p_2^\ell)^*\wt{L_2}=\wt{L}(\tau^{\ell }(\wt{\bd_2})-\wt{\bd_1},\frac{\wt{\lambda_2}^{r_2/g}}{\wt{\lambda_1}^{r_1/g}})
\] on $C^{hn}$. Then the dimension formula follows from Lemma \ref{lem: H0}.
\end{proof}

\begin{prop}\label{prop: end}
For $f\in \wt{S}^{k,+}_n$, denote by $V_f\in \cU^{+}_{L,k+1}$ a vector bundle that corresponds to $f$. Then
\[
\ell(f)=\dim(\End(V_f))-p(f).
\] 
\end{prop}

\begin{proof}
Note that $\ell(f)=\ell(fs_+)$.
We first prove the indecomposable case. In this case $p(f)=1$ and
\[
\mathfrak{C}_{fs_+}=\Big\{ \alpha, s_+^n\alpha,\ldots, s_+^{kn}\alpha \Big\}.
\]
It is clear that $f|_\alpha$ is monotone for all $\alpha\in \mathfrak{C}_f$.
Given $\alpha\neq \beta\in \mathfrak{C}_f$,  there exists $i\in \alpha$ and $j\in \beta$ such that
\[
i<j, f(j)<f(i)
\] if and only if 
\[
w_i^\beta-w_i^\alpha=w_{f(i)}^\beta-w_{f(i)}^\alpha=-1
\] and 
\[
s_{i+2,f(i)-1}(w^\beta-w^\alpha)>0.
\] Moreover 
\[
s_{i+2,f(i)-1}(w^\beta-w^\alpha)=\# \{j\in\beta| i<j, f(j)<f(i)\}.
\]
We set 
\[
N_{\beta,\alpha}=\sum_{i\in\Omega(\alpha)}  s_{i+2,f(i)-1}(w^\beta-w^\alpha),
\] where $\Omega(\alpha)=\{i<f(i)<f^2(i)<\ldots<f^{|\ol{\alpha}|}(i)=i+(k+1)n\}\subset \alpha$ is a fundamental domain under translation.
The infinite cyclic group generated by $s^n_+$ acts $\mathfrak{C}_f\times\mathfrak{C}_f\setminus \Delta$ with $k$ orbits represented by $(\alpha, s^{\ell n}_+\alpha)$ for $\ell=1,\ldots,k$. Then we have
\[
\ell(f)=\sum_{\ell=1}^k N_{s^{\ell n}_+\alpha,\alpha}.
\]
Now we apply Proposition \ref{prop: end} to the case when $V_1=V_2=V_f=B(\bd,\wt{\lambda})$. Then
\[
\dim(\End(V_f))= \sum_{\ell=1}^k \theta(\tau^\ell(\bd)-\bd, 1)+ 1. 
\]
Observe that 
\[
N_{s^{\ell n}_+\alpha,\alpha}=\theta(\tau^\ell(\bd)-\bd, 1).
\] So we prove that
\[
\dim(\End(V_f))=\ell(f)+1.
\]

Suppose $V=V_1\oplus\ldots\oplus V_p$ with $V_i$ indecomposable and non isomorphic. Let 
\[
\mathfrak{C}_{fs+}=\mathfrak{C}_1\sqcup \ldots\sqcup \mathfrak{C}_p
\]  be the orbit decomposition  for the action of the cyclic group $\{s_+^{nd}|d\in\Z\}$. Recall that
\[
\mathfrak{C}_i=\{\alpha_i, \ldots, s_+^{n\cdot(\omega(\alpha_i)-1)}\alpha_i\}
\] for $i=1,\ldots,p=p(f)$.
For $u\neq v\in\{1,\ldots,p\}$,
the diagonal action of $\{s_+^{nd}|d\in\Z\}$ on $\mathfrak{C}_u\times \mathfrak{C}_v$ has exactly $g_{u,v}$ orbits where $g_{u,v}=\gcd(\omega(\alpha_u), \omega(\alpha_v))$. If we fix $\alpha_u\in \mathfrak{C}_u$ and $\alpha_v\in \mathfrak{C}_v$ then the set of representatives is
\[
\{ (\alpha_u,s_+^{\ell n}\alpha_v)\Big| \ell=0,\ldots, g_{u,v}-1\}.
\] 
The number
\[
N_{v,u}:=\sum_{\ell=0}^{g_{u,v}-1}\sum_{i\in \Omega(\alpha_u)}  s_{i+2,f(i)-1}(w^{s_+^{\ell n}\alpha_v}-w^{\alpha_u})
\] counts the pairs $i<j$ such that $f(j)<f(i)$ with $i\in\alpha_u$ and $j\in\alpha_v$. Let $V_u=B(\bd^u,\wt{\lambda_u})$ and $V_v=B(\bd^v,\wt{\lambda_v})$ be the indecomposable summands that correspond to $\mathfrak{C}_u$ and $\mathfrak{C}_v$. Since $\bd^u\neq \tau^{\ell}(\bd^v)$ for $\ell=0, \ldots, g-1$,  by Proposition \ref{prop: end}
\[
N_{v,u}=\dim\Big(\Hom(V_u,V_v)\Big).
\] Applying the formula proved in the indecomposable case, we have
\[
\ell(f)=\sum_{u\neq v} N_{v,u}+ \Big(\sum_{\ell=1}^p \dim(\End(V_\ell))\Big)-p=\dim(\End(V))-p.
\]
\end{proof}

\section{Applications}\label{appl-sec}
\subsection{Applications to Poisson geometry of positroid varieties}
\begin{definition}(See \cite[Section 5]{KLS})\label{def: positroid}
Given $f\in B(k,n)$, let 
\[
X_f:=GL_k \Big\backslash \{M| f_M=f\}\subset G(k,n)
\] where $GL_k$ acts by row transformation. We call $X_f$ the \emph{(open) positroid variety} with affine permutation $f$.
\end{definition}

The first part of the next Theorem recovers the result of Goodearl and Yakimov \cite{GY} in the case of Grassmannians.

\begin{thm}\label{thm: main}
The stratification
\[
G(k,n)=\bigsqcup_{f\in B(k,n)} X_f
\] coincides with the stratification by $T$-leaves with respect to the standard Poisson structure $\pi_{d,V}$. Moreover, for each $f\in B(k,n)$ there is a surjective smooth morphism 
from $X_f$ to an algebraic torus of dimension $p(f)-1$ whose fibers are symplectic leaves of $\pi^\prime_{d,V}$. 
\end{thm}
\begin{proof}
Recall that there is an isomorphism between the coarse moduli space of $\cN_{L,k}$ and $G(k,n)$. Since $\cN_{L,k}$ is a $\G_m$ gerbe over a smooth projective scheme, the $0$-shifted Poisson structure descends to the Feigin-Odesskii Poisson structure on the coarse moduli space (see Section \ref{sec: FO}). By Theorem \ref{thm: extension}, the stack morphism $q: \cN_{L,k}\to \cU^{++}_{L,k+1}$ is surjective and $T$-equivariant. In Theorem \ref{thm: fibration} we have shown that the homotopy fiber of $q$ at a residue gerbe is connected and $0$-shifted symplectic. Then there is a bijection between $T$-leaves of $\cN_{L,k}$ and the $T$-orbits of $\cU^{++}_{L,k+1}$, which is bijective to $B(k,n)$ by Theorem \ref{thm: extension}. In Theorem \ref{thm: comparison} we prove that the FO bracket $2\Pi_{W,L}$ coincides with $\pi'_{d,V}$ which completes the proof of the first part. 

Fix $f\in B(k,n)$ and a vector bundle $V_f$ in the corresponding $T$-orbit. Denote by $v$ the point represented by $V_f$ in the vector bundle stack. 
By Theorem \ref{thm: fibration}, $X_f$ can be identified with the coarse moduli scheme of $\cN^v_{L,k}$, which admits a smooth surjective morphism
\[
q^v: \cN^v_{L,k}\to Tv
\] with connected and $0$-shifted symplectic fibers. Since the morphism $\cN^v_{L,k}\to \cN_{L,k}$ is a locally closed immersion, $\cN^v_{L,k}$ is a $\G_m$ gerbe over $X_f$. On the other hand, $Tv$ is a $T_1$-gerbe over an algebraic tori $T_2$ where $T_1$ can be identified with the stabilizer subgroup at $v$ (see proof of Theorem \ref{thm: fibration}). Therefore $q^v$ descends to a scheme morphism $X_f\to T_2$ that is smooth surjective whose fiber is connected and symplectic with respect to  $\Pi_{W,L}$. It remains to compute the dimension of $T_2$. Suppose $V_f=V_1\oplus\ldots\oplus V_{p(f)}$ is a splitting into sum of indecomposables. For an indecomposable bundle $B(\bd,\wt{\lambda})$ with $\bd=(d_1,\ldots,d_{rn})$ and $t\in T$,
\[
t^*B(\bd,\wt{\lambda})=B(\bd,\prod_{i=1}^n t_i^{-\sum_{j=0}^{r-1} d_{i+j n}}\wt{\lambda})
\] Since the $T$-action commutes with direct sum, $T_2$ has dimension $p(f)-1$. Here minus one comes from the determinant fixed condition.
\end{proof}

\begin{rmk} Theorem \ref{thm: main}
enhances the results of \cite{GY} by giving a description of how a $T$-leaf is fibered by symplectic leaves. The base of these symplectic fibrations are certain torus orbits of $\cU^{+}_{L,k+1}$ that glue together to the open substack $\cU^{++}_{L, k+1}$. Up to our limited knowledge, this is the first nontrivial example that the moduli stack of symplectic leaves of a partial flag variety is explicitly constructed.  It is an interesting question to ask if one can generalize our theorem to arbitrary partial flag varieties. This is open even for type $A$.
\end{rmk}

\begin{rmk}\label{rmk: image}
 The orbit closure-inclusion relation of the $T$-action on $\cU^{++}_{L,k+1}$ gives a geometric interpretation of the Bruhat order on $B(k,n)$.
\end{rmk}

\begin{cor}
The above theorems have several interesting consequences.
\begin{enumerate}
\item[$(1)$] $X_f$ is symplectic if and only if $p(f)=1$, i.e. $V_f$ is indecomposable;
\item[$(2)$] The cyclic permutation of branches (of $C$) acts on the set of positroid variteis as Poisson isomorphism and the Kummer involution that sends the $i$-th component to the $(n-i)$-th component, acts as an anti Poisson isomorphism.
\end{enumerate}
\end{cor}
\begin{proof}
$(1)$ is an immediate consequence of  Theorem \ref{thm: main}. By \cite[Theorem 7.5]{HP-Bos}, auto-equivalences, in particular, automorphisms of $C$ acts on the moduli stack of complexes and it scales the $0$-shifted Poisson structure. It is easy to check that in the case of cyclic permutation the scalar is $1$ and in the case of Kummer involution the scalar is $-1$.
\end{proof}

We finish the paper by giving some example of $T$-leaves for $G(2,4)$.
\begin{example}
Consider $n=4, k=2$ and the affine permutation $f: (1234)\mapsto (3456)$. The characteristic matrix $A_{fs_+}$ is 
\[
\left[
\begin{array}{cccccccccccccc}
& 1 & 0&0 & 1 &0 & 0 & 1 & 0 & 0 & 1 &0 & 0 &\\
\ldots & 0& 1& 0 & 0 & 1 & 0 & 0 & 1 & 0 & 0 & 1 & 0 & \ldots\\
& 0& 0 & 1 & 0 & 0 & 1 &0 & 0 & 1 & 0 & 0 & 1 & 
\end{array}
\right].
\]
In this example all three $\Z_f$ orbits are equivalent under the action of $s_+^4$, i.e. $p(f)=1$. Therefore the corresponding vector bundle $V_f$ is indecomposable. By Proposition \ref{prop: end},  
\[
\dim \End(V_f)= \ell(f)+p(f)=0+1=1,
\] i.e. $V_f$ is simple. As a consequence, $G(2,4)$ is generically symplectic under the twisted standard Poisson structure. The open symplectic leaf is $X_f$.

The affine permutation $f: (1234)\mapsto (5364)$ has characteristic matrix 
\[A_{fs_+}=
\left[
\begin{array}{cccccccccccccc}
& 1 & 0&1 & 1 &0 & 0 & 0 & 0 & \\
\ldots & 0& 0& 0 & 0 & 1 & 0 & 1& 1 & \ldots\\
& 0& 1 & 0 & 0 & 0 & 1 &0 & 0 & 
\end{array}
\right].
\] In this case, the first and second row are equivalent by $s_+^4$. The corresponding vector bundle splits into a direct sum of a rank 2 indecomposable bundle and a line bundle, i.e. $p(f)=2$. Again by Proposition \ref{prop: end}, 
\[
\dim \End(V_f)= \ell(f)+p(f)=3+2=5. 
\] In this case $X_f\cong \G_m$ and the symplectic leaves are single points.

The affine permutation $f: (1234)\mapsto (2349)$ has characteristic matrix 
\[A_{fs_+}=
\left[
\begin{array}{cccccccccccccc}
& 1 & 0&1& 0 &0 & 0 & 0 & 0 & \\
\ldots & 0& 0& 0 & 0 & 1 & 0 & 1 & 0 & \ldots\\
& 0& 1 & 0 & 1 & 0 & 1 &0 & 1 & 
\end{array}
\right].
\]
The corresponding vector bundle is the direct sum of a rank two indecomposable and a line bundle. One can show that 
\[
\dim \End(V_f)= \ell(f)+p(f)=3+2=5. 
\]  Since $f\not\in B(2,4)$, $V_f$ is not an extension of $\OO_C^2$ by $L$. This example shows the subtlety of the extension problem since we could not distinguish this case from the previous case by the functions $\ell(f)$ and $p(f)$.
\end{example}

\bibliographystyle{alpha}      
\bibliography{mybib}   

\end{document}